\documentclass[3p,times]{scrartcl}

\author{Alexander Herzog}

\usepackage{amssymb}
\usepackage{amsthm}
\usepackage{array}

\usepackage{colortbl}
\usepackage{color}
\usepackage{amsmath}
\usepackage{booktabs}
\usepackage{graphicx}
\usepackage{wrapfig}
\usepackage{algorithm}
\usepackage[noend]{algpseudocode}
\usepackage{authblk}
\usepackage[width=500pt,height=700pt]{geometry}

\definecolor{lightgray}{gray}{.85}

\newcommand{\New}{^\mathrm{(Newton)}}
\newcommand{\gMGF}{^\mathrm{(gMGF)}}

\newtheorem{defi}{Definition}
\newtheorem{thm}{Theorem}[section]
\newtheorem*{prf}{Proof}
\newtheorem{lem}{Lemma}[section]
\newtheorem{ex}{Example}

\makeatletter
\def\BState{\State\hskip-\ALG@thistlm}
\makeatother

\begin{document}

\title{An iterative equation solver with low sensitivity on the initial value}

\author{Alexander Herzog}

\affil{Munich University of Applied Sciences HM, Department of Engineering and Management, Germany}

\maketitle

\begin{abstract}
The objective of this publication is to reduce the sensitivity of iterative equation solvers on the initial value. To this end, at the hand of Newton's method, we exemplify how to reformulate the initial problem by means of a set of generalized moment generating functions. The approach allows to choose that very function, which is best approximated by a linear function and thus allows to set up an efficient iteration procedure. As a result of this, the number of iterations required to meet a given precision goal is significantly reduced in comparison to Newton's method especially for large deviations between the initial value and the actual root. At the hand of seven academic examples and three applications we demonstrate that the computing time of the discussed approach reveals a far lower susceptibility on the initial value when compared to results from Newton's method. This insensitivity offers the prospect to implement iterative equation solvers for applications with strict real-time requirements such as power system simulation or on-demand control algorithms on embedded systems with low computing power. We are confident that the devised methodology may be generalized to other well-established iteration algorithms.
\end{abstract}

\section{Introduction}
\label{Sec:Intro}
Solving nonlinear equations of the form
\begin{equation}\label{Problem}
	f(x)=0,
\end{equation}
where $f:\ I\subseteq\mathbb{R}\rightarrow\mathbb{R}$, has numerous applications in various fields of mathematics, natural sciences and engineering. Unfortunately, for an enormous amount of real-world problems, analytical treatments are not available, such that iterative algorithms yielding sufficiently accurate approximate solutions are of great importance. 

One prominent example for such a technique consists in Newton's method, where, starting from an initial guess $x_0$ close to the actual root $\zeta$, in the $n+1$-st iteration step an improved approximation 
\begin{equation}\label{Newton}
	x_{n+1}=x_{n}-\frac{f(x_{n})}{f'(x_{n})},\qquad n\in\mathbb{N}
\end{equation}
is calculated from the previous result together with the corresponding function value and its first derivative. It is well-known, that this algorithm is of convergence order $p=2$ for a simple root $\zeta$ with the error equation 
\begin{equation}\label{EENewton}
	e\New_{n+1}=\frac{f''(\zeta)}{2f'(\zeta)}\left(e\New_{n}\right)^2,
\end{equation}
where $e\New_n=x_n-\zeta$.

A multitude of methods have been suggested to improve this convergence order and the numerical efficiency $\eta$ respectively. The latter is defined as $\eta=p^{1/E},$ where $E$ marks the number of function evaluations per iteration~\cite{Ostro}. For one-point iterations the Kung-Traub conjecture establishes that the upper bound for $\eta$ is given by $p\leq2^{E-1}$~\cite{KungTraub}. Those methods which fulfill equality are called optimal.

One possibility to achieve an improvement of the convergence order consists in taking higher derivatives of the function $f(x)$ into account as e.g. is done within Halley's method or, more generally, Householder's method. In~\cite{ABBASBANDY2003887, BASTO2006468} the Adomian decomposition method was used to derive algorithms to solve Eq.~(\ref{Problem}) with cubic convergence. In spite of the higher convergence order of most approaches involving higher derivatives, they are rarely used in practice due to the additional computing time as this often renders these algorithms numerically inefficient. One alternative in using the Adomian decomposition method to derive a cubically converging iteration but solely involving the function $f$ and its first derivative was suggested by Chun in~\cite{CHUN2006415}. As a matter of fact, this algorithm is a special case of iterations which can be obtained from the homotopy analysis method~\cite{Chun2006}. The latter allows to derive iterations with convergence order three and four. Interestingly, the homotopy analysis method allows to solve Eq.~(\ref{Problem}) by Newton-like iterations in arbitrarily many ways. 

A different approach to tackle Eq.~(\ref{Problem}) by iteration consists in multi-step algorithms such as Traub's method~\cite{Traub}. Further, more recent examples for such schemes were discussed in~\cite{9462087}, as well as in~\cite{OZBAN2022100157}, where a whole family of fourth-order methods has been proposed. In~\cite{THANGKHENPAU2023100243} two derivative-free optimal families of bi-parametric iterations were suggested and studied at the hand of engineering applications. Starting from Newton's method in~\cite{MCDOUGALL201420} an algorithm was suggested where the argument for the derivative is determined by a predictor-corrector approach to increase the convergence order to $p=1+\sqrt{2}$.

In general, also for the two-steps methods approximating higher order derivatives is a quite common strategy. Different opportunities have been discussed in the literature. Here, we like to mention~\cite{MAHESHWARI2009383} where the second order derivative was approximated by a combination of $f$ and $f'$ leading to an iteration with $p=4$. The two-step method suggested in~\cite{9709846} is kept free of derivatives and also displays fourth-order convergence.  In~\cite{math7010055} [1,n]-order Pad\'e approximations were used to systematically obtain iterations involving higher derivatives. To circumvent the actual calculation of the latter, approximations were used thus changing the procedures to two-step methods.

 An approach to increase the convergence order without the necessity to treat higher derivatives, consists in taking higher order terms of $f(x_n)$ and $f'(x_n)$ into account as was suggested in~\cite{math7070655} at the hand of eight methods. Systematic treatments in the same direction utilize weight functions, which are chosen such that next-to-leading order expressions within the investigated iteration vanish. Such strategies are broadly used both for one-point methods, see e.g.~\cite{CHUN20071103, axioms8020055}, as well as for multi-step methods with and without memory as discussed in~\cite{SHARMA2023100270, Sharma, DZUNIC20124917, ZHANLAV2017414, WANG2018710, math8010108, sym15010228, doi:10.1080/27690911.2022.2130914, sym14102020}. 
 
A drawback common to most of the iterations available is that, albeit they are of a high local convergence order, their performance strongly depends on whether the initial guess $x_0$ is close to the sought for root $\zeta$. For some applications this is however not tolerable. For instance, the real-time simulation of power systems is an important tool to ensure their faultless operation. Here, typical targets for the computing time lie in the sub-millisecond range~\cite{9576157}. Due to the rapid fluctuations of the electromagnetic fields, in such situations, a convenient choice of the initial value is hard to determine in advance. Further examples are given by a number of algorithms implemented on automotive electronic control units designated for on-board diagnosis or control purposes which have to be executed in real time. Here, the solution of the equation under consideration can vary from one time step to another by several orders of magnitude, see c.f.~\cite{en13051148, Perrone_2022, Vagapov2020}. Also under these circumstances it may by no means be clear, what a proper choice for an initial value might be. If the performance of the chosen iteration sensitively depends on this, it therefore is not ensured that the value calculated from $x_0$ results in convergence of the iteration while fulfilling the real-time requirement. To make matters worse, the embedded systems used in the applications mentioned above very often only exhibit very reduced computing powers. On the other hand, for these applications the demands for the precision of an iteration is comparatively low. However, especially safety-relevant functions are to be executed dependably within a short amount of time. From this standpoint, a fast convergence irrespective of the initial guess is of high interest for these type of problems.

Therefore, in this contribution, we intend to discuss a strategy for iterative approaches to approximately solve Eq.~(\ref{Problem}), which does not aim so much at increasing the local convergence order. Instead, our objective is to decrease the overall number of iterations necessary to reach convergence, thus reducing the sensitivity of the computing time on the initial value. Thus, we try to focus on convergence for the whole range which is to be considered for a given problem. Such iterative techniques may serve as enablers for rapid power system simulations or on-demand control algorithms on embedded systems with low computing power as pursued in the automotive industry or process engineering. It is our belief, that the suggested algorithm can be used to generalize most of the iterations mentioned in this introduction. However, for the sake of simplicity, we discuss the procedure with the aim to generalize Newton's method, such that this famous algorithm serves as the benchmark for our approach throughout the publication.

Our method is remotely related to the aforementioned systematical derivation from the homotopy analysis method discussed by Chun~\cite{Chun2006} in that we could view the approach as a reformulation of Eq.~(\ref{Problem}). There is also some connection to the approaches discussed in~\cite{math7070655} and our procedure can be viewed as a generalization of the methods four and five mentioned in this reference.

Our procedure is based on the reformulation of Eq.~(\ref{Problem}) in such a way, that it can be represented by a formal moment generating function (MGF), as was suggested for control problems in~\cite{Vagapov2020}. A subsequent generalization then allows to define infinitely many generalized MGFs (gMGFs), which has been suggested again for the solution of control problems, c.f.~\cite{Herzog2020}. Here, we adapt this approach to Eq.~(\ref{Problem}). Having defined the gMGFs, we choose the one representation, which exhibits the most rapid convergence behavior, if approximated by a linear function. From this, a procedure requiring a small number of iterations for the calculation of the approximate root can be deduced. Thus, our method allows to bring the considered function to a form, which can be regarded as approximately linear. In many cases, this allows to reduce the number of iterations considerably. This scheme comes at the cost of evaluating several exponential or logarithmic functions, which at a first glance seems to compromise the efficiency of the method. However, as is well-known, a number of algorithms exist, which allow to represent these functions to high accuracy at very low computational costs~\cite{schraudolph1999fast, cawley2000fast, muller2020elementary, moroz2022simple, PERINI201837}. Thus, in contrary to existing approaches, which aim at capturing higher orders through higher derivatives or multi-step methods, the main idea of our approach is to cast Eq.~(\ref{Problem}) into a form which allows for an essential linear approximation which converges rapidly. To the best of our knowledge, such a methodology has not been devised before.

Since the discussion of $p$ and $\eta$ is usually restricted to the leading order, which may however be inappropriate if $x_0$ is far away from $\zeta$, we discuss our method at the hand of the gain in distance which we obtain in processing the iteration. In order to compare our procedure with Newton's method we investigate the respective residuals and gains in distance as a function of the computing time connected to this for both algorithms Moreover, we study the computational order of convergence (COC). As we shall see at the hand of seven academic examples and three real-world applications, for initial values far away from the sought-after root, the gMGF approach exhibits computing times which turn out to be up to an order of magnitude lower than within Newton's method.

The paper is organized as follows: In Section~\ref{Const}, we construct the iterative gMGF approach on the basis of Newton's method. A detailed analysis on the convergence order together with a comparison to Newton's iteration and improvements for applied implementation is discussed in Section~\ref{Analysis}. We shall also give a step by step algorithm sketch by means of pseudo-code complemented by a flowchart visualizing the approach. Numerical examples are given in Section~\ref{NumEx}. Here, we will compare the number of iterations, the number of function evaluation, as well as the computing time of our approach with Newton's method. Moreover, some detailed insights into the convergence behavior by means of the residuals shall be addressed. In Section~\ref{Apps}, we compare the performance of our model with the one obtained from Newton's method at the hand of three industrial applications where the real-time behavior is of importance. Section~\ref{SumOut} summarizes the results of our investigations and gives an outlook on applications and further research activities.

\section{Construction of the iterative generalized moment generating function method}\label{Const}
In this section, we intend to develop our algorithm. To this end we use the following definitions:
\begin{defi}
The gMGF of degree 0 reads
\begin{equation}\label{H0}
	\mathcal{H}_0(x_n,\xi_n)=-\sigma\left(f(x_n)\right)\left(f(x_n+\xi_n)-f(x_n)\right),
\end{equation}
where $\sigma(z)$ gives the sign of the real number $z$ and $\xi_n=x_{n+1}-x_n$ with $n\in\mathbb{N}\cup\{0\}$.
\end{defi}
\begin{defi}
We define the gMGF of arbitrary degree $k\in\mathbb{Z}$ recursively via
\begin{equation}\label{Hk}
	\mathcal{H}_k(x_n,\xi_n)=\ln\left[1+\mathcal{H}_{k+1}(x_n,\xi_n)\right].
\end{equation}
\end{defi}
From a Taylor expansion to order $N$ with respect to $\xi_n$ we then obtain
\begin{equation}\label{HkTaylor}
	\mathcal{H}_k(x_n,\xi_n)=\sum_{l=1}^Nh_{k,l}(x_n)\frac{\xi_n^l}{l!}+\mathcal{O}\left(\xi_n^{N+1}\right),
\end{equation}
where we have introduced the generalized moments of degree $k$ and order $l$ via
\begin{equation}\label{Moments}
	h_{k,l}(x_n)=\left.\frac{\partial^l\mathcal{H}_k(x_n,\xi_n)}{\partial\xi_n^l}\right|_{\xi_n=0}.
\end{equation}

Formally the power series given in Eq.~(\ref{HkTaylor}) resembles the form of a MGF known from statistics. There, the moments allow to obtain characterizing quantities of the probability density function such as the mean value, the variance or the skewness. A number of statistical problems simplify when represented by the corresponding cumulant generating function (CGF) which is defined as the natural logarithm of the MGF. Thus, it becomes apparent from Eq.~(\ref{Hk}) that $\mathcal{H}_k$ can be interpreted as a formal modified CGF of $\mathcal{H}_{k+1}$. It is moreover well-known that cumulants and moments are interrelated to each other, which therefore also applies in the present context. In particular, the moments can be expressed by the cumulants through the equation of Faà di Bruno, and a similar relation should hold for the moments of the gMGFs with adjacent degrees. In fact, from the recursion given in Eq.~(\ref{Hk}) we may establish

\begin{thm}\label{FaaBrunoGen}
	The generalized moments $h_{k,l}(x_n)$ of $\mathcal{H}_k(x_n,\xi_n)$ and the generalized moments $h_{k+1,l}(x_n)$ of $\mathcal{H}_{k+1}(x_n,\xi_n)$ obey the relations
	\begin{equation}\label{MomnMomnm1}
		h_{k,l}(x_n)=\sum_{m=1}^lB_{l,m}\left(h_{k-1,l}(x_n),h_{k-1,l}(x_n),\dots,h_{k-1,l}(x_n)\right),
	\end{equation}
	and
	\begin{equation}\label{Momnm1Momn}
		h_{k-1,l}(x_n)=\sum_{m=1}^l(-1)^{m-1}(m-1)!B_{l,m}\left(h_{k,l}(x_n),h_{k,l}(x_n),\dots,h_{k,l}(x_n)\right),
	\end{equation}
	where the functions $B_{l,m}$ represent the incomplete Bell polynomials.
\end{thm}
\begin{prf}
	From Eq.~(\ref{Hk}), we observe that $\mathcal{H}_k(x_n,\xi_n)$ resembles the form of the CGF with respect to the formal MGF $\mathcal{H}_{k+1}(x_n,\xi_n)$, where the only difference in comparison to the conventional relation between the two lies in the argument of the logarithm which in our case is given by $1+\mathcal{H}_{k+1}(x_n,\xi_n)$ instead of $\mathcal{H}_{k+1}(x_n,\xi_n)$. This discrepancy however does not play a role for the calculation of the respective generalized moments (c.f. Eq.~(\ref{Moments})) as can be verified by direct differentiation. Therefore, the generalized moments of $\mathcal{H}_{k}(x_n,\xi_n)$ and $\mathcal{H}_{k+1}(x_n,\xi_n)$ fulfill the same relations among each other as the moments calculated from a MGF and the cumulants obtained from the CGF. These relations are given by Eqs.~(\ref{MomnMomnm1}) and~(\ref{Momnm1Momn}) as was proven in \cite{Comtet}.
	\qed
\end{prf}
Theorem~\ref{FaaBrunoGen} allows us to establish relations of the generalized moments of arbitrary degree to those of degree zero. For the present purpose we solely discuss the relations for the moments of the first two orders in the following.
\begin{lem}\label{lemh1}
	For the first order generalized moments of degree 0 and degree $k$ the following holds:
	\begin{equation}\label{h1}
		h_{k,1}(x_n)=h_{0,1}(x_n).
	\end{equation}
\end{lem}
\begin{prf}
	 We have $B_{1,1}(x)=x$. Inserting this into Eq.~(\ref{MomnMomnm1}) yields $h_{k,1}(x_n)=h_{k-1,1}(x_n)$. Eq.~(\ref{h1}) follows from recursion.
	\qed
\end{prf}
\begin{lem}\label{lemh2}
	Between the second order generalized moment of degree $k$ and the second order generalized moment of degree $0$ the following holds:
	\begin{equation}\label{h2}
		h_{k,2}(x_n)=h_{0,2}(x_n)+kh_{0,1}^2(x_n).
	\end{equation}
\end{lem}
\begin{prf}
	 For $k=0$ the statement contained in Eq.~(\ref{h2}) is trivial. In the case, that $k>0$ and setting $l=2$ in Eq.~(\ref{MomnMomnm1}) we obtain 
	 \begin{equation}
	 	h_{k,2}(x_n)=h_{k-1,2}(x_n)+h_{k-1,1}^2(x_n).
	\end{equation}
	Moreover, from Eq.~(\ref{h1}), we readily conclude that
	 \begin{equation}\label{hk2Lem}
	 	h_{k,2}(x_n)=h_{k-1,2}(x_n)+h_{0,1}^2(x_n).
	 \end{equation}
	 Then an expression analogous to Eq.~(\ref{hk2Lem}) applies for $h_{k-1,2}(x_n)$. Thus, inserting $h_{k-1,2}(x_n)$, $h_{k-2,2}(x_n),\dots, h_{k-q,2}$ consecutively, we end up at
	 \begin{equation}
	 	h_{k,2}(x_n)=h_{k-q,2}(x_n)+qh_{0,1}(x_n),
	 \end{equation}
	 which, letting $q=k$ proofs Eq.~(\ref{h2}) for $k>0$. The proof for $k<0$ can be carried out in full analogy.
	\qed
\end{prf}
With the gMGFs (c.f. Eq.~(\ref{Hk})), we utilize Eq.~(\ref{HkTaylor}) where we choose the degree $k$ such that the curvature of $\mathcal{H}_{k}(x_n,\xi_n)|_{\xi_n=0}$ is minimal. In the vicinity of $x_n$ this leads us to the gMGF which may be represented best by a linear function. This optimal degree $\varkappa$ amounts to the gMGF for which the absolute value of the second order moment is minimal. From Eq.~(\ref{h2}) we readily obtain
\begin{equation}\label{kappa}
	\varkappa=-r\left(\frac{h_{0,2}(x_n)}{h_{0,1}^2(x_n)}\right),
\end{equation}
where $r(a)$ rounds $a$ to the respective closest integer value.

Having calculated $\varkappa$, we use the power series in Eq.~(\ref{HkTaylor}) of optimal degree $\varkappa$ to linear order to obtain
\begin{equation}\label{LinRep}
	\mathcal{H}_\varkappa(x_n,\xi_n)=h_{\varkappa,1}(x_n)\xi_n+\mathcal{O}\left(\xi_n^2\right).
\end{equation}
We now substitute the left-hand side of this equation with the value $\mathfrak{H}_\varkappa(x_n)$ which we obtain by setting $f(x_n+\xi_n)=0$ in the definition of $\mathcal{H}_\varkappa(x_n,\xi_n)$, i.e.
\begin{equation}\label{DefVal}
	\mathfrak{H}_\varkappa(x_n)=\left.\mathcal{H}_\varkappa(x_n,\xi_n)\right|_{f(x_n+\xi_n)=0}.
\end{equation}
From this definition, we obtain
\begin{equation}
	\mathfrak{H}_0(x_n)=\left|f(x_n)\right|
\end{equation}
for degree zero and the recursion
\begin{equation}
	\mathfrak{H}_k(x_n)=\ln\left[1+\mathfrak{H}_{k+1}(x_n)\right]
\end{equation}
applies. Inserting Eq.~(\ref{DefVal}) into Eq.~(\ref{LinRep}) we obtain
\begin{equation}\label{Hakappa}
	\mathfrak{H}_\varkappa(x_n)=h_{\varkappa,1}(x_n)\xi_n+\mathcal{O}\left(\xi_n^2\right),
\end{equation}
which, neglecting higher order terms and utilizing Eq.~(\ref{h1}), allows us to calculate $x_{n+1}$ as
\begin{equation}\label{xjp1General}
	x_{n+1}=x_n+\frac{\mathfrak{H}_\varkappa(x_n)}{h_{0,1}(x_n)},
\end{equation}
which is the iteration equation to the method proposed in this manuscript.

\section{Analysis and improvements for implementation}\label{Analysis}
A closer inspection of the algorithm proposed in this publication reveals, that for $\varkappa=0$ the iteration given in Eq.~(\ref{xjp1General}) is identical to Newton's method. For $\varkappa\neq0$ it may however be regarded as a generalization to this iteration. We also like to point out, that for $\varkappa=1$ the suggested approach is equivalent to method four in~\cite{math7070655} for single roots, whereas method 5 of the same reference is captured by $\varkappa=-1$ also for the case of single roots.

From this, the convergence order of the proposed algorithm should equal the one of Newton's method, which is proven in the following theorem.
\begin{thm}
	The method proposed in Eq.~(\ref{xjp1General}) with Eq.~(\ref{kappa}) is of quadratic convergence. If $x_n$ is sufficiently close to the actual root $\zeta$, then the error $e\gMGF_{n+1}$ of iteration step $n+1$ defined via  
	\begin{equation}\label{Errjp1}
		e\gMGF_{n+1}=x_{n+1}-\zeta
	\end{equation}
	fulfills the error equation
	\begin{equation}\label{errorEq}
		e\gMGF_{n+1}=\frac{h_{\varkappa,2}(\zeta)}{2h_{0,1}(\zeta)}\left(e\gMGF_n\right)^2.
	\end{equation}
\end{thm}
\begin{prf}
	Inserting Eq.~(\ref{xjp1General}) into Eq.~(\ref{Errjp1}) we obtain
	\begin{equation}\label{ejp1andej}
		e\gMGF_{n+1}=x_n+\frac{\mathfrak{H}_\varkappa(x_n)}{h_{0,1}(x_n)}-\zeta.
	\end{equation}
	From a Taylor expansion a straightforward calculation yields
	\begin{equation}\label{TaylormathfrakH}
		\mathfrak{H}_\varkappa(x_n)=-h_{0,1}(\zeta)e\gMGF_n-\frac{1}{2}h_{\varkappa,2}(\zeta)\left(e\gMGF_n\right)^2+\mathcal{O}\left(\left(e\gMGF_n\right)^3\right).
	\end{equation}
	Analogously, we expand $h_{0,1}(x_n)$ in Eq.~(\ref{ejp1andej}) to linear order as
	\begin{equation}\label{Taylorh01}
		h_{0,1}(x_n)=h_{0,1}(\zeta)+h_{0,2}(\zeta)e\gMGF_n+\mathcal{O}\left(\left(e\gMGF_n\right)^2\right).
	\end{equation}
	Inserting Eq.~(\ref{TaylormathfrakH}) and Eq.~(\ref{Taylorh01}) into Eq.~(\ref{ejp1andej}) and using $e\gMGF_n=x_n-\zeta$ we readily end up at
	\begin{equation}
		e\gMGF_{n+1}=e\gMGF_n-\frac{e\gMGF_n+\frac{1}{2}\frac{h_{\varkappa,2}(\zeta)}{h_{0,1}(\zeta)}\left(e\gMGF_n\right)^2+\mathcal{O}\left(\left(e\gMGF_n\right)^3\right)}{1+\frac{h_{0,2}(\zeta)}{h_{0,1}(\zeta)}e\gMGF_n+\mathcal{O}\left(\left(e\gMGF_n\right)^2\right)}.
	\end{equation}
	Evaluating this expression to lowest order, we obtain
	\begin{equation}
		e\gMGF_{n+1}=\frac{h_{\varkappa,2}(\zeta)}{2h_{0,1}(\zeta)}\left(e\gMGF_n\right)^2,
	\end{equation}
	which completes the proof.
	\qed
\end{prf}
For $\varkappa=0$, Eq.~(\ref{xjp1General}) is identical to the iteration formula given by Eq.~(\ref{Newton}). Consequently the well-known convergence behavior of Newton's method is captured by Eq~(\ref{errorEq}) in this case.

We infer that from the above analysis the proposed method is only marginally superior to Newton's method, in the sense, that the coefficient relating $e_n$ to $e_{n+1}$ is smaller than for the proposed method if $\varkappa\neq0$. Combined with the fact that our method requires a number of additional function evaluations to determine $\varkappa$ as well as for the calculation of $\mathfrak{H}_\varkappa(x_n)$, at first sight, Newton's method appears to overall be superior to it.

Under some conditions however the proposed method may indeed offer benefits. One drawback of Newton's method rests in the fact that its performance strongly depends on the starting value $x_0$ for the iteration. This may turn out critical if the initial guess $x_0$ for the root of $f(x)$ is far away from $\zeta$. Under these circumstances the number of needed iterations within the gMGF algorithm may be considerably smaller than within Newton's method, as the nonlinear manipulation of the function evaluation induced by Eq.~(\ref{Hakappa}) favors a bigger step-size $|x_{n+1}-x_n|$ during the evaluation. If, due to this effect, the number of iterations within Newton's method exceeds the ones used within the proposed algorithm drastically, the computing time of our approach may be considerably smaller than the one achieved by the iteration given in Eq.~(\ref{Newton}), as will be shown in Sections.~\ref{NumEx} and~\ref{Apps}.

From an application standpoint for a given problem the method with the better performance should be chosen. This is however fulfilled automatically within the proposed method, because, as explained above, for $\varkappa=0$ the gMGF algorithm coincides with Newton's method. A drawback of our method is that the calculation of the degree can be rather time consuming. To avoid the comparatively costly explicit calculation of $\varkappa$ from Eq.~(\ref{kappa}) during the iteration of the algorithm, we may evaluate $|\varkappa|\in\mathbb{N}\cup\{0\}$ prior to its implementation and deploy the chosen degree by means of staggered if-then-else constructs. Together with the determination of the sign of $f(x_n)$ we thus obtain the degree more efficiently as compared to the calculation suggested in Eq.~(\ref{kappa}). A closer analysis reveals that this scheme reduces the computing time of the method considerably. Thus, the algorithm suggested in this publication is finalized. An outline is given in \textbf{Algorithm~\ref{Algoritm1}}.

\begin{figure}
	\includegraphics[scale=0.55]{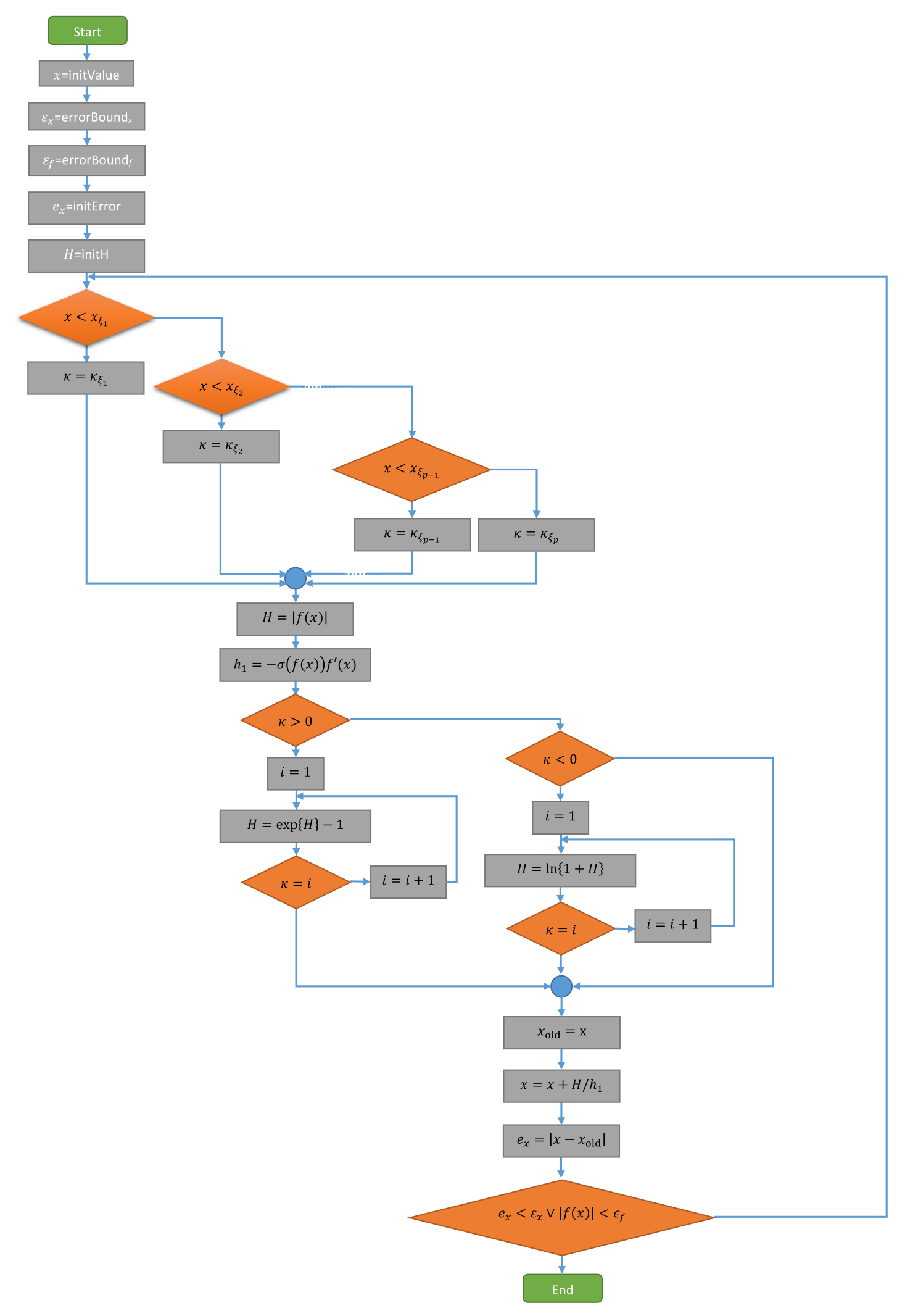}
	\caption{Flowchart visualizing the generalized moment generating function approach: After the initialization (first five processes) $\varkappa$ is determined in the staggered if-else construct. After this evaluation the generalized moment generating function of degree zero and the first moment are calculated. Depending on the sign of $\varkappa$ the corresponding generalized moment generating function is evaluated. After storing the result from the previous step, the corresponding updated value as well as the error is evaluated in order to decide, whether the iteration continues or not.}\label{Flowchart}
\end{figure}

Moreover, in Fig.~\ref{Flowchart} the approach is visualized by means of a flowchart: After the initialization (first five processes) $\varkappa$ is determined in the staggered if-else construct. After this evaluation $\mathfrak{H}_0$ and the first moment are calculated. Depending on the sign of $\varkappa$ the corresponding gMGF is evaluated. After storing the the result from the previous step, the corresponding updated value as well as the error is evaluated in order to decide, whether the iteration continues or not.

\begin{algorithm}[t]
\caption{Iterative generalized moment generating function approach}\label{Algoritm1}
\begin{algorithmic}[1]
\State $x\leftarrow$ initValue, $\epsilon_x\leftarrow$ errorBound$_x$, $\epsilon_f\leftarrow$ errorBound$_f$, $e_x\leftarrow$ initError, $H\leftarrow$ initH 
\While {$e_x>\epsilon_x$ \textbf{or} $H>\epsilon_f$}

\textit{Comment: Calculate the optimal degree:}
\If {$x<x_\mathrm{\xi_1}$}
	$\varkappa\leftarrow\varkappa_\mathrm{\xi_1}$
	\Else
	\If {$x<x_\mathrm{\xi_2}$}
		$\varkappa\leftarrow\varkappa_\mathrm{\xi_2}$		
	\Else
	\If{ $\dots$}\dots
	\State\vdots
	\If{ $x<x_\mathrm{\xi_{p-1}}$}
		$\varkappa\leftarrow\varkappa_\mathrm{\xi_{p-1}}$
	\Else\ 
	$\varkappa\leftarrow\varkappa_\mathrm{\xi_p}$		
	\EndIf
	\EndIf
\EndIf
\EndIf

\textit{Comment: Calculate H$_0$, and the first moment:}
\State $H\leftarrow|f(x)|$, h1$\leftarrow -\sigma(f(x))f'(x)$

\textit{Comment: Calculate H$_\varkappa$; if $\varkappa=0$ no operations have to be performed.}
\If {$\varkappa>0$}
	\For {degIter=$1,2,\dots,\varkappa$}
		\State H$\leftarrow\exp\{\textnormal{H}\}-1$
	\EndFor
	\Else
	\If {$\varkappa<0$}
		\For {degIter=$-1,-2,\dots,\varkappa$}
		\State H$\leftarrow\ln\{1+\textnormal{H}\}$
	\EndFor
	\EndIf
\EndIf
\State xOld$\leftarrow$x

\textit{Comment: Iteration equation}
\State x$\leftarrow$x+H/h1
\State $e_x=|\textnormal{x}-\textnormal{xOld}|$
\EndWhile
\end{algorithmic}
\end{algorithm}

Addressing the possible drawbacks of the proposed method, as already stated, the calculation of $\mathfrak{H}_\varkappa$ necessitates additional function evaluations of exponential or logarithmic functions. However, we like to indicate, that a number of very efficient methods exist to this end (c.f~\cite{schraudolph1999fast, cawley2000fast, muller2020elementary, moroz2022simple, PERINI201837}), such that the additional computational cost connected to this is compensated for sufficiently big values of $|e_0^{(\gamma)}|$ with $\gamma\in\{\mathrm{Newton},\mathrm{gMGF}\}$.

In addition, a possible limitation of the method may occur for $\varkappa>0$ and $\mathfrak{H}_0\geq1$, as in this case the exponential term stemming from the inversion of Eq.~(\ref{Hakappa}) comes drastically into play and turns out to be the more severe the higher the chosen degree is. This circumstance may necessitate to limit $\varkappa$ to an upper value appropriate for the problem at hand. Under these circumstances the full performance of the gMGF method can not be utilized. However, in such a case, the truncated value for the degree may still require less computation time as compared to Newton's method. We like to stress, however, that we did not observe a problem of this kind during the evaluation of the numerical experiments discussed in Sections~\ref{NumEx} and~\ref{Apps}. 

Finally, if the problem stated in Eq~(\ref{Problem}) is to be solved for a periodic function, in the case of an increased step-size, the iteration may result in a root which was not the expected one. This problem may also occur within Newton's method. However, due to the bigger step sizes taken for $\varkappa\neq0$ as a trend, the effect is more likely to occur within the gMGF approach. To avoid such a behavior, the degree would have to be truncated again, thus inhibiting the method to develop its full power. 

To evaluate the performance of the present approach in comparison to Newton's method, the following definition is given:
\begin{defi}
	The gain $g^{(\gamma)}_n$ of iteration step $n$ is defined by the reduction of the magnitude of the error $|e^{(\gamma)}_n|$ for the $n$-th iteration step in comparison to its predecessor. Hence,
	\begin{equation}
		g^{(\gamma)}_n=\left|e^{(\gamma)}_{n-1}\right|-\left|e^{(\gamma)}_{n}\right|.
	\end{equation}
\end{defi}

\section{Numerical examples}\label{NumEx}
In this section, we will demonstrate the performance of the iterative gMGF method by addressing the root finding to several functions. In order to compare its performance to Newton's method, instead of varying the initial value, we choose to vary the dependent variable, i.e. we intend to calculate the inverse function $x=f^{-1}(y)$ to the equation
\begin{equation}\label{InversFunc}
	y=f(x),
\end{equation}
where $f(x)$ is chosen from categories of nonlinear equations which regularly occur in numerical calculations such as exponential functions, trigonometric functions, power functions, polynomials or rational functions and combinations thereof. We point out, that varying $y$ while keeping $x_0$ fixed amounts to a variation of the distance between the initial value and the root $\zeta$.
Subsequently, $\delta y$ defines the step size of the function value $y$. Both methods will operate with the same initial values $x_0$ respectively. For the given range $y\in[y_\text{min},y_\text{max}]$ for both methods $x=f^{-1}(y)$ will be illustrated. The iteration terminates, once $|x_n-x_{n-1}|,|y-f(x_n)|\leq10^{-15}$ is reached. This comparatively crude value is chosen with a focus on applications on embedded systems with low computing power, like process control computers or automotive electronic control units. Moreover the number of iterations, $I$, as well as the number of function evaluations, $E$, until convergence each as a function of $y$ shall be displayed. In the framework of the gMGF method, the latter is given by
\begin{equation}
	E\gMGF=\sum_{n=0}^{I-1}\left(3+\left|r(\varkappa(x_n)\right|\right).
\end{equation}
For all of the discussed examples, we will address the computing time $t(y)$ needed for convergence defined as
\begin{equation}
	t(y)=\sum_{n=1}^{I(y)}\tau_n,
\end{equation}
where $\tau_n$ is the computing time of iteration step $n$. In order to minimize computational fluctuations of the processing device (Intel(R) Core(TM) i5-8265U CPU @ 1.60GHz   1.80 GHz with an installed RAM of 16,0 GB and a 64-Bit-based processor), $t(y)$ is determined from the mean value of 1000 independent iteration procedures for each method. For all iterations the language Wolfram Mathematica 13.2.0.0 was used.

In the subsequent graphical representations, the black thick curve depicts results for Newton's method, whereas the red thin line displays the corresponding quantities for the algorithm introduced in this publication. The function $x=f^{-1}(y)$ is shown in subplot (a), while (b) and (c) respectively represent $I$ and $E$. Finally subplot (d) shows the computing time $t(y)$.

Moreover, for every example a comparison of the convergence behavior of Newton's iteration and our method is given for one value of $y$. To this end, the magnitude of the error $|e_n|$ as well as the gain $g_n$ as a function of the respective computing time are displayed for one starting point for both methods. These respective quantities are displayed in the subplots (e) and (f), where the black circles and red diamonds show results from Newton's iteration and the iterative gMGF-method each. In addition, for a given value of $y$ we visualize the quantity $|x_n-x_{n-1}|$ alongside the COC $\varrho_n$ defined by
\begin{equation}
	\varrho_n=\frac{\log\left|\frac{x_n-\zeta}{x_{n-1}-\zeta}\right|}{\log\left|\frac{x_{n-1}-\zeta}{x_{n-2}-\zeta}\right|}
\end{equation}
for all iteration steps until the aforementioned precision goal is reached. This information will be displayed in tables.

Subsequently, we discuss the individual examples:
\begin{ex}\label{Exam1}
	We investigate the function
	\begin{equation}
		y=x e^{x^2}-\sin^2x+3\cos x+5,
	\end{equation}
	where we consider $y\in[-10,8]$ with $\delta y=0.1$. The initial value for the iteration is chosen as $x_0=0$.
\end{ex}
\begin{figure}[t]
	\begin{center}
		\includegraphics[width=.95\linewidth]{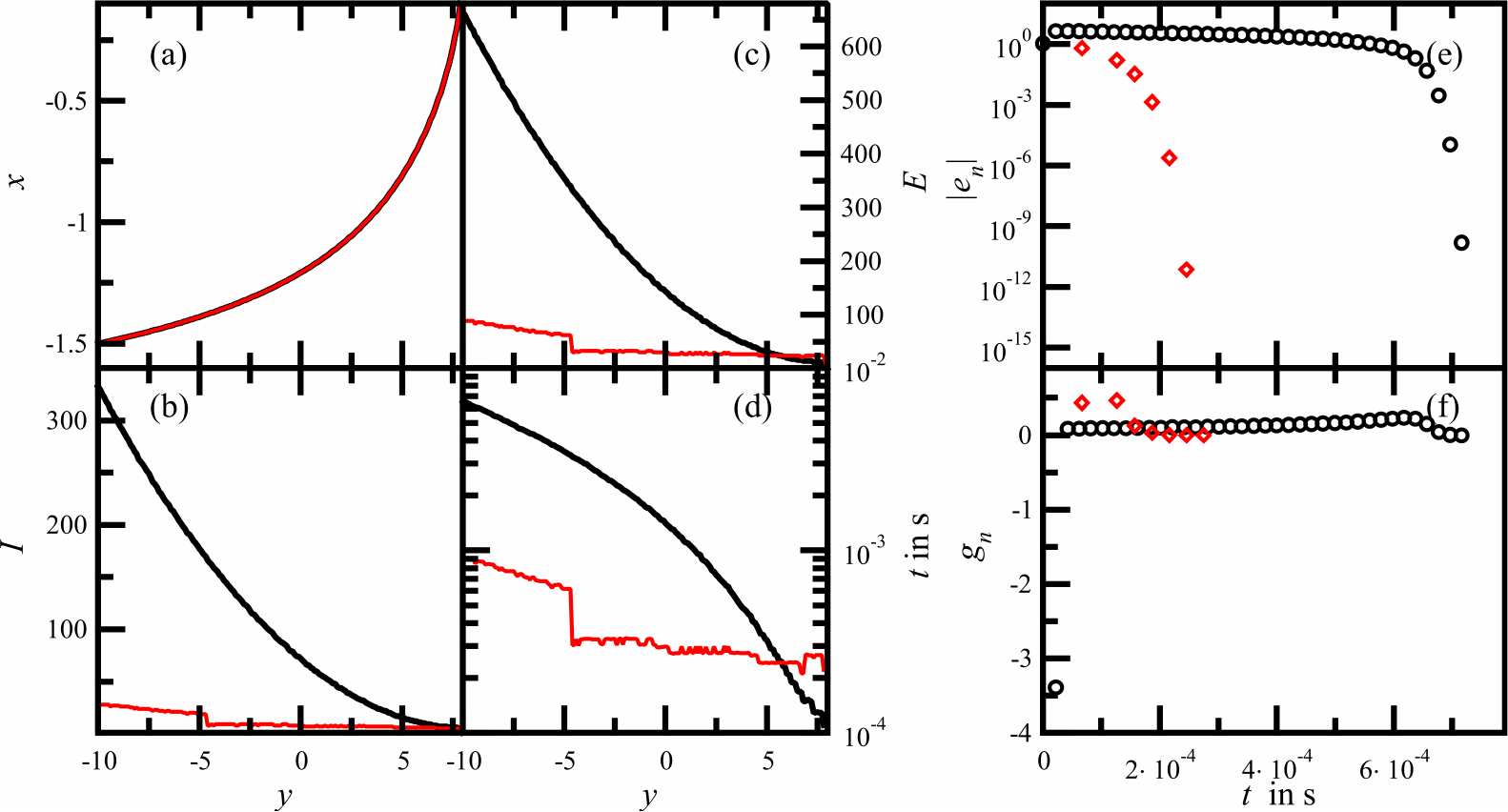}
		\caption{Comparison of results for \textbf{Example 1} as obtained from Newton's method (black curve) and our approach (red curve): subplot (a) displays the inverted function $x=x(y)$, whereas (b) shows the number of iterations needed for convergence. The number of function evaluations connected to this can be seen from subplot (c). In subplot (d) the computing time until convergence is shown. Especially when the starting value deviates strongly from the actual solution our approach outperforms Newton's method. Subplots (e) and (f) show criteria to compare the convergence behavior of Newton's method and our approach for \textbf{Example 1} with $y=2.5$. The circles display these quantities as obtained from Newton's method, whereas the diamonds are calculated from the approach suggested in this publication: In subplot (e) the magnitude of the residuals as a function of the computing time is shown. The respective gains are displayed in subplot (f). While Newton's method displays an increase of the error within the first iteration step, the high gain in the early stage of the iteration is responsible for the rapid convergence of our method.}\label{Ex1}
	\end{center}
\end{figure}
Fig.~\ref{Ex1} reveals the results obtained for \textbf{Example}~\ref{Exam1}. We observe that both iterations capture the correct inverse function (c.f. subplot (a)). For the complete interval of $y$ taken into account, we find $I\gMGF\leq I\New$ as can be seen from subplot (b). Moreover, we see, that $I\gMGF$ remains comparatively moderate with a maximum value of 28 within the considered parameter range, while $I\New$ increases considerably faster with increasing distance between $x_0$ and the actual value $x=f^{-1}(y)$.  The highest value observed in that case is located at $y=10$ and amounts to $I\New=332$.

To the contrary within the domain $[5.6,8]$ the number of functional evaluations $E^\mathrm{(Newton)}=2I^\mathrm{(Newton)}$ (c.f. subplot (c)) fulfills $E^\mathrm{(Newton)}\leq E\gMGF$. This is due to the fact, that obtaining $\varkappa(x_{n-1})$ from the if-else-construct contributes to $E^\mathrm{(gMGF)}$ on the one hand. On the other hand,  calculating $\mathfrak{H}_\varkappa(x_{n-1})$ for $\varkappa\ne0$ involves additional functional evaluations. Still for $y\leq5.5$, $I^\mathrm{(Newton)}$ exceeds  $I^\mathrm{(gMGF)}$ to an extent, that the relation $E^\mathrm{(Newton)}>E^\mathrm{(gMGF)}$ holds. The maximum for $E^\mathrm{(gMGF)}=89$ in the given interval occurs at $y=-10$, whereas Newton's method exhibits a maximum value which is over seven times higher.

The computing time displayed in subplot (d) shows a similar behavior: For $y=8$, Newton's method outperforms our approach by a factor of two. This difference however gradually decreases with decreasing $x$ and vanishes at $y=5.5$. Below this value, our algorithm displays a lower computing time. As for the number of iterations and the number of functional evaluations, $t^\mathrm{(gMGF)}\leq0.91437$ ms remains moderate with the maximum value occurring at $y=-10$. The maximum computing time for Newton's method is given by $t^\mathrm{(Newton)}(y=-10)\approx6.5972$ ms. Hence, as the number of function evaluations, the computing time for Newton's method is over seven times higher than for our approach at this point. A closer investigation revealed, that the difference $\Delta t(y)=t\New(y)-t\gMGF(y)$ between the computing times of the two methods grows quadratic and persists beyond the investigated interval.

Subplot~(e) gives some insight into the convergence behavior of Newton's method and the gMGF approach respectively for $y=2.5$. Interestingly, we observe that within the first iteration, the magnitude of the residual resulting from Newton's method increases compared to the initial residual, i.e. $|e_1^\mathrm{(Newton)}|>|e_0^\mathrm{(Newton)}|$. For the subsequent iteration steps we then find $|e_{n+1}\New|<|e_n\New|$. However, the circumstance, that $|e_1^\mathrm{(Newton)}|>|e_0^\mathrm{(Newton)}|$, seems to be responsible for the large number of iterations, which, in this case is given by $I\New(y=2.5)=36$. We conclude, that even though the initial value is sufficiently close to the actual solution for the iteration to converge, it is a bad choice in terms of computing time. From this perspective, the sensitivity of Newton's method with respect to the starting value is displayed in this example. As can be seen, no such difficulty does occur in the present case for the iterative gMGF approach. For the latter, we observe that $|e_{n+1}\gMGF|< |e_n\gMGF|$, such that the algorithm converges to the precision goal within seven iterations with $e_7\gMGF=8.83\cdot10^{-18}$. We point out, that the single iterations demand a higher computing time. However, due to the rapid convergence with each iteration step, the gMGF algorithm outperforms Newton's method for the given example. We also observe, that the difference in the overall computing time tends more to the favor of the gMGF approach, the lower the precision goal is chosen. From a practical viewpoint, this underlines that the suggested algorithm is particularly suitable for implementation in applications, where a dependable computing time over a large variable range is of importance, whereas the precision goal is comparatively low. This is for instance realized for diagnosis or control  in the automotive or process engineering industry. The corresponding algorithms are often implemented on embedded systems of low computation power, with a demand for real-time response. On the other hand, typical quantities addressed in this area are temperature,  pressure, enthalpy, current or voltage. It is quite rare, that these quantities have to be known with more than five significant digits. The same can be said about power supply systems, which often also require rapid computing with a comparatively low precision goal.

Finally in subplot~(f),  the gain for the respective iteration steps is displayed. Within Newton's method the negative gain, i.e. the loss, during the first iteration step is clearly visible. For $n>1$ we observe $g_n\New\geq0$ with a maximum value of $g_{30}\New\approx0.23601$. To the contrary we have $g_n\gMGF>0$ for all iteration steps taken into account, thus hinting towards the circumstance, that the initial value is ''properly'' chosen for the used approach. The maximum value is given by $g_2\gMGF\approx0.46415$. Both, the fact, that this almost exceeds $g_{30}\New$ by a factor of two as well as the circumstance that the maximum occurs at an early iteration step, can be viewed as evidence that for the investigated value, the iterative gMGF method is superior to Newton's iteration.

In Tab.~\ref{tabular:T_MSL1}, another view on the convergence behavior of the two methods taken into consideration is shown for $y=7$. We find that Newton's method requires 7 iterations to reach convergence whereas $I\gMGF(7)=5$ holds. Moreover, we see, that with respect to the number of iterations the gMGF approach shows a more rapid convergence. Also, we find, that $\varrho_n\New$ starts with a negative value. This again hints towards the fact, that within the first iteration step, Newton's method yields a deterioration with respect to the initial value. For $n\geq3$ we observe that the COC tends to 2. The corresponding expression is far less regular within the gMGF algorithm: The rapid convergence at the early stage of the iteration, yields $\varrho_2\gMGF=3.02$ followed by a decrease after which the value predicted in Eq~(\ref{errorEq}) occurs. Finally, we observe that due to the similar numbers of iterations, the computing time within Newton's method is lower than within the gMGF approach for the studied case.

\begin{table}[t]
	\centering
	\caption{Numerical results for~\textbf{Example 1:} $y=x e^{x^2}-\sin^2x+3\cos x+5$ with $y=7$ and $x_0=0$}
	\label{tabular:T_MSL1}
	\begin{tabular}[!ht]{c c c c c c c c c c c c c c c}
	\toprule[2pt]
		Newton & $n=1$ & $n=2$ & $n=3$ & $n=4$ & $n=5$ & $n=6$ & $n=7$ & $I$ & CPU time\\
		\midrule
		$|x_n-x_{n-1}|$ & $1$ & $0.33$ & $0.19$ & $4.2\cdot10^{-2}$ & $1.9\cdot 10^{-3}$ & $3.7\cdot10^{-6}$ & $ 1.4\cdot10^{-11}$ & $7$ & $1.53\cdot 10^{-4}$ s\\
		$\varrho_n$ &  & $-3.64$ & $1.88$ & $1.90$ & $1.98$ & $2.00$ & $2.00$\\
		\midrule\midrule
		gMGF & $n=1$ & $n=2$ & $n=3$ & $n=4$ & $n=5$ &  &  & $I$ & CPU time\\
		\midrule
		$|x_n-x_{n-1}|$ & $0.30$ & $0.13$ & $4.1\cdot10^{-3}$ & $1.4\cdot10^{-5}$ & $1.7\cdot 10^{-10}$ &   &   & $5$ & $2.68\cdot 10^{-4}$ s\\
		$\varrho_n$ &  & $3.02$ & $1.62$ & $1.99$ & $2.00$\\
		\bottomrule[2pt]
	\end{tabular} 
\end{table} 

\begin{ex}\label{Exam2}
	This example is given by
	\begin{equation}
		y=x e^{x^2}-\sin^2x+3\cos x+5,
	\end{equation}
	with $y\in[8.1,100]$ with the step size $\delta y=0.1$. As the initial value we use $x_0=1$.
\end{ex}
\begin{figure}[t]
	\begin{center}
		\includegraphics[width=.99\linewidth]{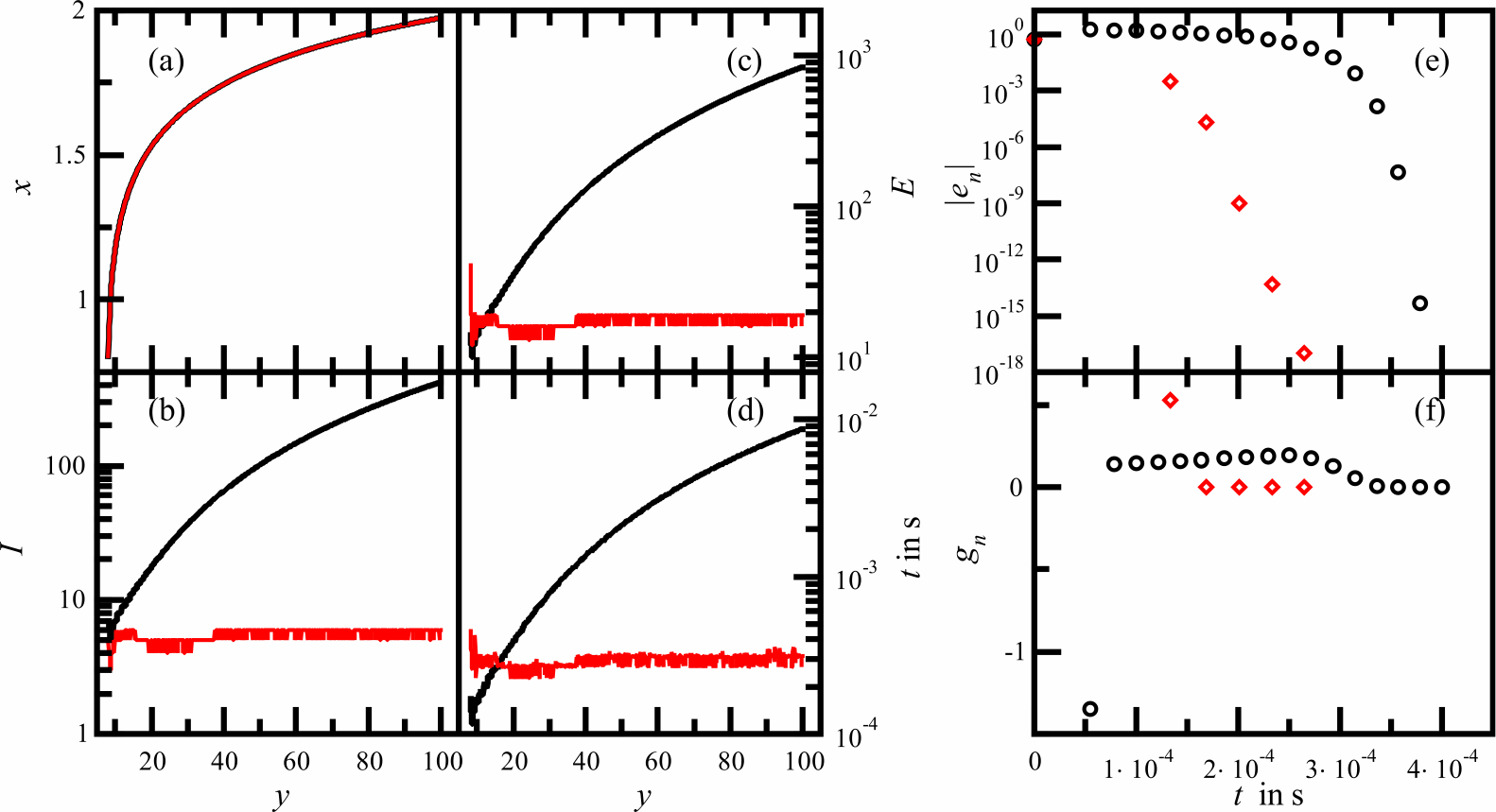}
		\caption{Comparison of results for \textbf{Example 2} as obtained from Newton's method (black curve) and our approach (red curve): subplot (a) displays the inverted function $x=x(y)$, whereas (b) shows the number of iterations needed for convergence. The number of function evaluations connected to this can be seen from subplot (c). In subplot (d) the computing time until convergence is shown. Especially when the starting value deviates strongly from the actual solution, our approach outperforms Newton's method. Subplots (e) and (f) show criteria to compare the convergence behavior of Newton's method and our approach for \textbf{Example 2} with $y=20$. The circles display these quantities as obtained from Newton's method, whereas the diamonds are calculated from the approach suggested in this publication: In subplot (e) the residual as a function of the computing time is shown. The respective gains are displayed in subplot (f). While Newton's method displays an increase of the error within the first iteration step, the high gain in the early stage of the iteration is responsible for the rapid convergence of our method.}\label{Ex2}
	\end{center}
\end{figure}
\begin{table}[t]
	\centering
	\caption{Numerical results for~\textbf{Example 2:} $y=x e^{x^2}-\sin^2x+3\cos x+5$ with $y=20$ and $x_0=1$.}
	\label{tabular:T_MSL}
	\begin{tabular}[!ht,width=.98\linewidth]{c c c c c c c c c c c c c c c}
	\toprule[2pt]
		Newton & $n=1$ & $n=2$ & $n=3$ & $n=4$ & $n=5$ & $n=6$ & $n=7$ \\
		\midrule
		$|x_n-x_{n-1}|$ & $2.41$ & $0.14$ & $0.15$ & $0.15$ & $0.16$ & $0.17$ & $ 0.18$ \\
		$\varrho_n$ &  & $-6.14\cdot10^{-2}$ & $1.13$ & $1.14$ & $1.17$ & $1.19$ & $1.23$\\\midrule
		&  $n=8$ & $n=9$ & $n=10$ & $n=11$ & $n=12$ & $n=13$ & $n=14$ & &\\\midrule
		$|x_n-x_{n-1}|$ & $0.18$ & $0.19$ & $0.19$ & $0.18$ & $0.13$ & $5.5\cdot10^{-2}$ & $8.0\cdot10^{-3}$ & \\
		$\varrho_n$ & $1.28$ & $1.34$ & $1.43$ & $1.55$ & $1.70$ & $1.85$ & $1.96$\\\midrule
		&  $n=15$ & $n=16$ & $n=17$ &  &    &  $I$ & CPU time &\\\midrule
		$|x_n-x_{n-1}|$ & $1.4\cdot10^{-4}$ & $4.6\cdot10^{-8}$ & $4.4\cdot10^{-15}$ &  &  & $17$ &   $4.00\cdot 10^{-4}$ s  & \\
		$\varrho_n$ & 2.00 & $2.01$ & $2.00$ & \\\midrule\midrule
		gMGF & $n=1$ & $n=2$ & $n=3$ & $n=4$ & $n=5$  & $I$ & CPU time\\
		\midrule
		$|x_n-x_{n-1}|$\ \  & $0.53$ & $3.1\cdot10^{-3}$ & $2.2\cdot10^{-5}$ & $1.0\cdot10^{-9}$ & $2.2\cdot 10^{-13}$   & $5$ & $2.65\cdot 10^{-4}$ s\\
		$\varrho_n$ &  & $0.971$ & $2.00$ & $1.52$ & $2.00$\\
		\bottomrule[2pt]
	\end{tabular} 
\end{table}

Results for \textbf{Example~\ref{Exam2}} can be found in Fig.~\ref{Ex2}. Both methods capture the inverse function to high accuracy as can be observed from subplot~(a). The number of iterations, displayed in subplot~(b) reveals that the number of iterations is smaller for our method for all values of $y$ taken into consideration. Here, within the gMGF approach $3\leq I\gMGF\leq6$ holds, while we observe that the corresponding quantity strongly increases for Newton's method for $y>10$. Overall, we have $I\New\in[5,420]$, where the minimum value occurs for $y\in[8.5,8.8]$.

Correspondingly, the number of function evaluations within Newton's method ranges between 10 and 840 (c.f. subplot~(c)). To the contrary we find $E\gMGF(y=8.1)=42$ corresponding to $I\gMGF=5$. Apart from this and the value $y=8.2$, where $E\gMGF=24$ in the considered range, we have $E\gMGF\leq19$.

The difference in the number of function evaluations is also displayed for the computing times depicted in subplot~(d). However, we observe that from this perspective, Newton's method can be regarded as superior in a broader range for $y$ as compared to the function evaluations. We find $t\New(y)<t\gMGF(y)$ for $y\in[8.1,16.1]$, whereas our approach outperforms Newton's method regarding the computing time for $y\geq16.2$. Summarizing, we find, that $t\New(y)$ increases disproportionate with increasing $|y(\zeta)-y(x_0)|$ where for the considered range we observed the maximum value $t\New(100)\approx8.01389$ ms. Contrary to this, the computing time needed for the gMGF algorithm to converge to the precision is less sensitive with respect to $y$, where the maximum computation time amounts to $t\gMGF(8.1)=0.45212$ ms. As in the previous example, $\Delta t(y)$ increases quadratically for high values of $y$.

In subplot~(e) we investigate the convergence behavior for $y=20$. As in the previous example, we observe, that $|e_1\New|>|e_0\New|$. During the subsequent iterations, the magnitude of the errors decrease, i.e. $|e_{n+1}\New|<|e_n\New|$ for $n>2$. In total 17 iterations are needed to reach the precision goal within Newton's method. For the gMGF method, we observe a strictly decreasing absolute value of the residuals, such that in total four iterations are required for convergence. For the first iteration, we detect a high computing time, which is due to the evaluation of $\mathfrak{H}_\varkappa(1)$. Again, we observe, that a lower precision goal favors the gMGF approach as far as computing time is concerned.

From subplot (f), we again observe the loss of Newton's method within the first iteration step. Apart from this, we have $g_n\New>0$ with a maximum value of $g_{10}\New\approx0.19269$. Contrary to this, the maximum gain for the gMGF method occurs within the first iteration reading $g\gMGF_1=0.52979$. For the considered case however this first iteration obviously takes a particularly long time which almost corresponds to five iterations within Newton's method. Still, for the considered value, we have $t\gMGF<t\New$.

Finally, in Tab.~\ref{tabular:T_MSL} further information on the convergence of the two methods is shown for $y=20$. As already discussed at the hand of Fig~\ref{Ex2}, Newton's method displays a loss in the first iteration step from which a negative COC for $n=2$ results. After this, a slight increase in $|x_n-x_{n-1}|$ is observed for $n\leq10$ which is also displayed in the increase of $\varrho\New_n$. For $n>10$, we see, that $|x_n-x_{n-1}|$ decreases, while the COC tends towards its final value, 2.00. Conversely, in the framework of the gMGF, we again observe, that convergence is very rapid already in the first iteration steps. Again, the COC within this approach does not show a clear trend. It reaches 2.00 in the final iteration step.
\begin{ex}\label{Exam4}
	Consider the function rule
	\begin{equation}
		y=x^{1/3}\left(x- e^{x}\right).
	\end{equation}
	Let us investigate the interval $y\in[-80,-0.5],$ with the step size $\delta y=0.1$ and the initial value $x_0=0.5$.
\end{ex}
\begin{figure}[t]
	\begin{center}
		\includegraphics[width=.99\linewidth]{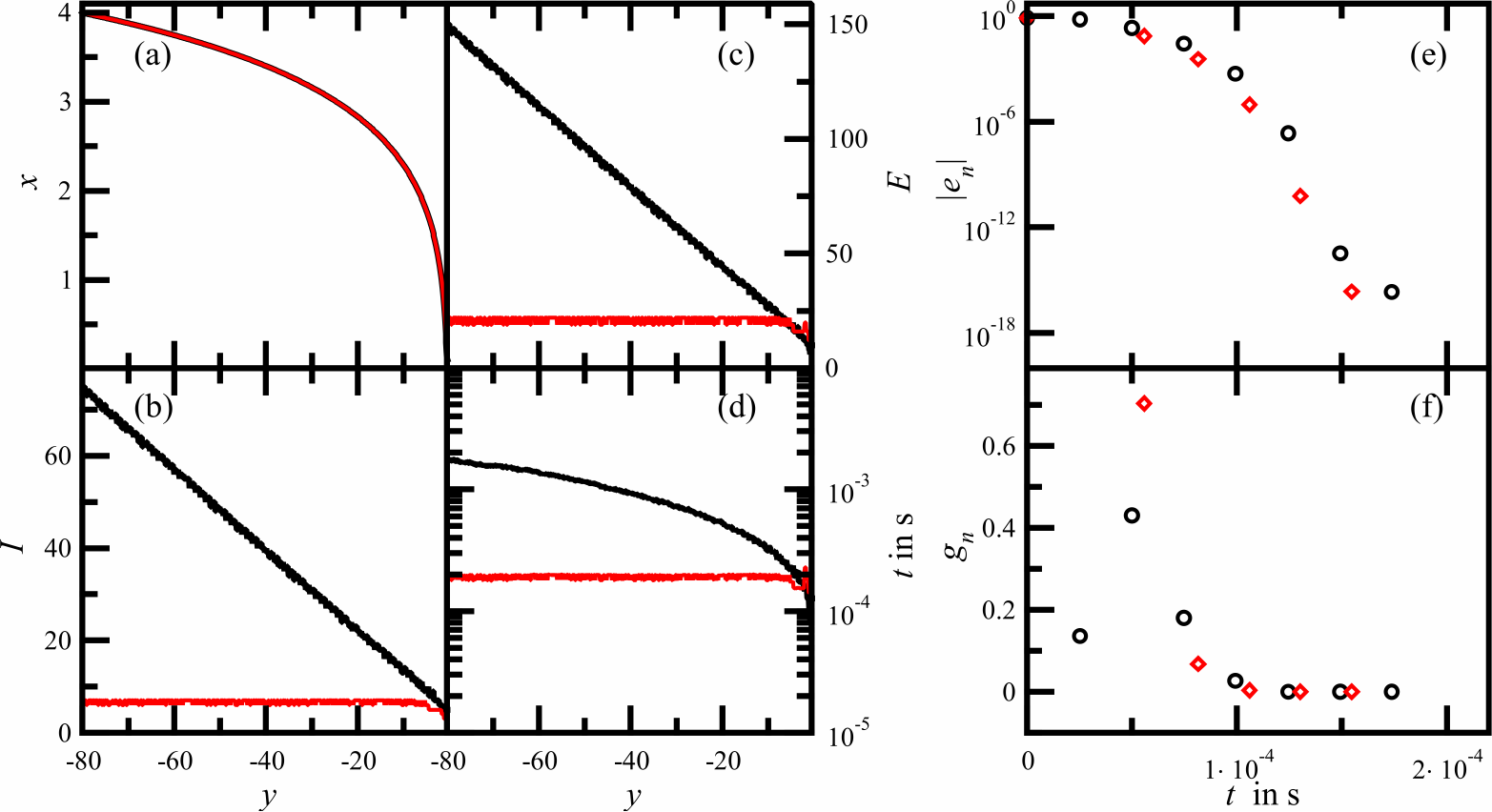}
		\caption{Comparison of results for \textbf{Example~\ref{Exam4}} as obtained from Newton's method (black curve) and our approach (red curve): subplot (a) displays the inverted function $x=x(y)$, whereas (b) shows the number of iterations needed for convergence. The number of function evaluations connected to this can be seen from subplot (c). In subplot (d) the computing time until convergence is shown. The gMGF approach exhibits an essentially constant computing time. Especially when the starting value deviates strongly from the actual solution, our approach outperforms Newton's method. Subplots (e) and (f) show criteria  to compare the convergence behavior of Newton's method and our approach for \textbf{Example~\ref{Exam4}} with $y=-2.5$. The circles display these quantities as obtained from Newton's method, whereas the diamonds are calculated from the approach suggested in this publication: In subplot (e) the residual as a function of the overall computing time is shown. The respective gain of the methods is depicted in subplot (f).}\label{Ex4}
	\end{center}
\end{figure}
\begin{table}[t]
	\centering
	\caption{Numerical results for~\textbf{Example 3:} $y=x^{1/3}(x- e^{x})$ with $y=-10$ and $x_0=0.5$}
	\label{tabular:T_MSL4}
	\begin{tabular}[!ht]{c c c c c c c c c c c c c c c}
	\toprule[2pt]
		Newton & $n=1$ & $n=2$ & $n=3$ & $n=4$ & $n=5$ & $n=6$&  $n=7$ \\
		\midrule
		$|x_n-x_{n-1}|$ & $8.1$ & $0.96$ & $0.95$ & $0.94$ & $0.92$ & $0.86$& $0.76$ \\
		$\varrho_n$ &  & $-0.131$ & $1.19$ & $1.23$ & $1.28$ & $1.35$& $1.46$  \\\midrule
		&$n=8$& $n=9$ & $n=10$ & $n=11$ &$n=12$ &$n=13$ & $I$ & CPU time & &\\\midrule
		$|x_n-x_{n-1}|$ & $0.57$& $0.29$ & $6.0\cdot10^{-2}$ & $2.3\cdot10^{-3}$ & $3.3\cdot10^{-6}$ &$6.5\cdot10^{-13}$& $13$&$3.14\cdot10^{-4}$ s\\
		$\varrho_n$ & $1.60$ & $1.77$ & $1.93$ & $1.99$ & $2.00$ & $2.00$ \\\midrule\midrule
		gMGF & $n=1$ & $n=2$ & $n=3$ & $n=4$ & $n=5$ & $n=6$ & $I$ & CPU time\\
		\midrule
		$|x_n-x_{n-1}|$\ \  & $2.1$ & $0.23$ & $3.8\cdot10^{-2}$ & $9.0\cdot10^{-4}$ & $4.9\cdot 10^{-7}$ & $1.4\cdot10^{-13}$ & $6$ & $1.79\cdot 10^{-4}$ s\\
		$\varrho_n$ &  & $1.02$ & $1.95$ & $1.99$ & $2.00$ & $2.00$\\
		\bottomrule[2pt]
	\end{tabular} 
\end{table}
\begin{table}[t]
	\centering
	\caption{Numerical results for~\textbf{Example 4:} $y=\frac{1}{x^2}+\frac{10}{x^4}+\frac{100}{x^{10}}$ with $y=1.5$ and $x_0=2.5$}
	\label{tabular:T_MSL5}
	\begin{tabular}[!ht]{c c c c c c c c c c c c c c c}
	\toprule[2pt]
		Newton & $n=1$ & $n=2$ & $n=3$ & $n=4$ & $n=5$ & $n=6$&  $n=7$ &$n=8$\\
		\midrule
		$|x_n-x_{n-1}|$ & $1.9$ & $0.065$ & $0.072$ & $0.079$ & $0.088$ & $0.099$& $0.11$ & $0.13$\\
		$\varrho_n$ &  & $-0.110$ & $1.18$ & $1.19$ & $1.22$ & $1.25$& $1.30$ & $1.38$ \\\midrule
		& $n=9$ & $n=10$ & $n=11$ &$n=12$ &$n=13$ & $n=14$&$n=15$&$n=16$\\\midrule
		$|x_n-x_{n-1}|$ & $0.14$ & $0.16$ & $0.14$ & $0.071$ &$0.013$& $3.4\cdot10^{-4}$&$2.3\cdot10^{-7}$&$1.0\cdot10^{-13}$\\
		$\varrho_n$  & $1.49$ & $1.64$ & $1.80$ & $1.92$ & $1.98$ & $2.00$&$2.00$&$2.00$\\\midrule
		& $n=17$& & & & & &  $I$ & CPU time\\\midrule
		$|x_n-x_{n-1}|$ & $2.2\cdot10^{-16}$ &  &  &  & & & $17$&$3.91\cdot10^{-4}$ s\\
		$\varrho_n$  & $2.00$ \\\midrule\midrule
		gMGF & $n=1$ & $n=2$ & $n=3$ & $n=4$ & &  & $I$ & CPU time\\
		\midrule
		$|x_n-x_{n-1}|$\ \  & $0.75$ & $0.065$ & $9.2\cdot10^{-4}$ & $3.3\cdot10^{-8}$ &  & & $4$ & $2.83\cdot 10^{-4}$ s\\
		$\varrho_n$ &  & $1.50$ & $2.47$ & $1.98$ & \\
		\bottomrule[2pt]
	\end{tabular} 
\end{table}

Results for \textbf{Example~\ref{Exam4}} are depicted in Fig.~\ref{Ex4}. As in the previous examples, both methods yield the inverse function for the whole range considered as can be seen from subplot~(a). Subplot~(b) shows, that $I\gMGF(y)\leq I\New(y)$ for all considered values of $y$. We find that $I\New$ increases approximately linear with increasing $|x_0-\zeta|$ with a maximum value of 75, whereas $I\gMGF\in[3,7]$. The same tendencies are displayed for the number of function evaluations shown in subplot (c). However, due to additional evaluations to determine $\varkappa$ and to calculate $\mathfrak{H}_\varkappa$ for $\varkappa\neq0$, $E\New(y)\leq E\gMGF(y)$ for $y>-3.4$, whereas within the rest of the considered interval, we have $E\gMGF(y)<E\New(y)$.

Subplot~(d) also reveals an approximately linear increase of the computing time for increasing $|x_0-\zeta|$ for Newton's method. The maximum value amounts to $1.8543$ ms. To the contrary, within the iterative gMGF approach, the computing time remains essentially constant with $t\gMGF\approx0.2$ ms over the considered parameter range. As far as computing time is concerned, the gMGF approach outperforms Newton's method for $y\leq-2.3$. Hence, even though we have $E\New(y)\leq E\gMGF(y)$ for $y>-3.4$, we still find $t\New(y)>t\gMGF(y)$ for $y\in[-80.0,-2.3]$. This may be explained by the efficient evaluation of $\varkappa$ by means of the if-else construct which is connected to very low computing times. A closer investigation of the computing times, revealed that $\Delta t(y)$ grows linearly with $y$. We have checked, that this tendency persists beyond the  interval considered in Fig.~\ref{Ex4}.

To investigate the convergence behavior by means of the residuals and the gain, we use $y=-2.5$. For this value, the gMGF algorithm performs slightly better than Newton's method as far as the computing time is concerned. The number of iterations respectively read $I\New(-2.5)=7$ and $I\gMGF(-2.5)=5$. From subplot (e), we observe that the magnitudes of the residuals decrease for every iteration step, within both methods. Moreover, we find, that the first step of the gMGF method takes slightly more than the computing time necessary for two evaluations within Newton's iteration. Still, the improvement obtained from this step is quite remarkable, as can be seen from the ratio of the initial residual and the corresponding successor which amounts to over an order of magnitude, i.e. $e_0/e_1\gMGF\approx10.8523$.

Considering the gain from subplot~(f), we find that $g_1\New(-2.5)+g_2\New(-2.5)<g_1\gMGF(-2.5)$. Interestingly, even though the gain within Newton's method for the considered example is never negative, we still have $g_1\New(-2.5)<g_2\New(-2.5)$. A closer analysis of the obtained data reveals, that for $y<-2.9$ the gain of the first iteration step within Newton's method turns negative, thus revealing the sensitivity of the algorithm on the initial value also for \textbf{Example~\ref{Ex4}}. To the contrary, the gMGF approach does not show such a sensitivity for the considered problem, such that the computing time remains approximately the same over the studied range.

Additional information on the convergence behavior is contained in Tab.~\ref{tabular:T_MSL4} for $y=-10$. From $\varrho_2\New<0$, we may deduce, that the first iteration step within Newton's method is connected to a massive loss, the magnitude of which can be estimated from the expression $|x_1-x_0|$ shown in the table. For $n\geq2$, we see, that $|x_n-x_{n-1}|$ strictly decreases and it takes 13 iterations until convergence is reached. As far as the COC is concerned, we observe that it strictly increases and after the initial negative value, it converges to the theoretical value 2.00 during the iterations. Conversely $\varrho_n\gMGF>0$ within the six iterations which are required for convergence and no loss is detected within this approach. Instead we again find a convergence of the COC to its theoretical value with a faster rate, than observed within Newton's method. From the upper row of the gMGF-quantities listed we see that in particular the first iteration step causes considerable progress. A closer investigation shows, that the singular iterations within the gMGF algorithm require more time, than those of Newton's method. However, the difference in the number of iterations leads to a smaller overall CPU time for our suggested approach.
\begin{ex}\label{Exam5}
As the fourth example, we intend to invert the function
	\begin{equation}
		y=\frac{1}{x^2}+\frac{10}{x^4}+\frac{100}{x^{10}},
	\end{equation}
where we consider $y\in[0.1,10000]$ with the step size $\delta y=0.01$. As initial value, we choose $x_0=2.5$.
\end{ex}
\begin{figure}[t]
	\begin{center}
		\includegraphics[width=.99\linewidth]{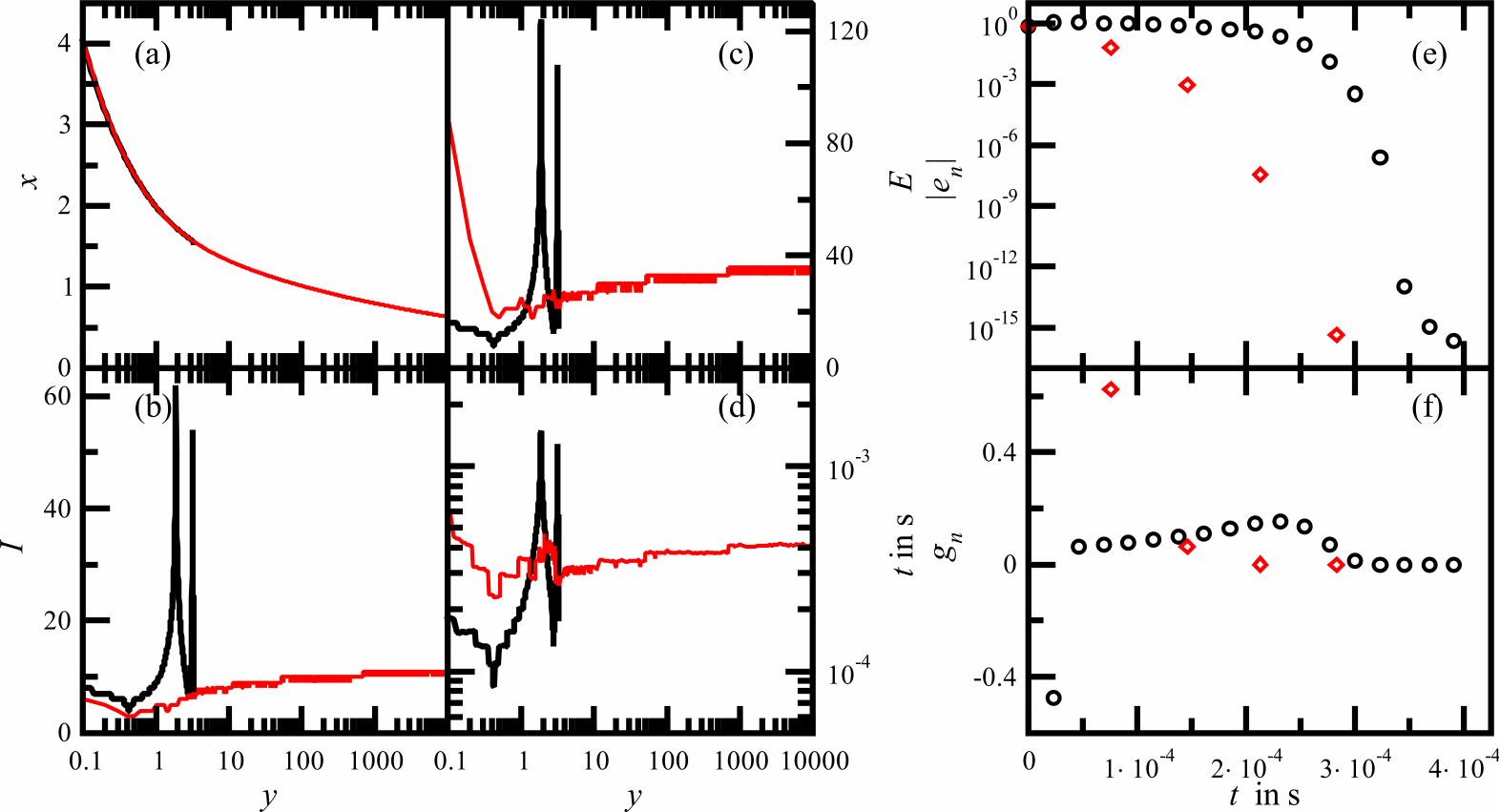}
		\caption{Comparison of results for \textbf{Example~\ref{Exam5}} as obtained from Newton's method (black curve) and our approach (red curve): subplot (a) displays the inverted function $x=x(y)$. We find, that for $y\geq 3.33$, Newton's method does diverge, whereas no such problem is observed for the gMGF method. Subplot (b) shows the number of iterations needed for convergence. The number of function evaluations connected to this can be seen from subplot (c). In subplot (d) the computing time until convergence is shown. Subplots (e) and (f) show criteria to compare the convergence behavior of Newton's method and our approach for \textbf{Example~\ref{Exam5}} with $y=1.5$. The circles display these quantities as obtained from Newton's method, whereas the diamonds are calculated from the approach suggested in this publication: In subplot (e) the residual as a function of the respective computing time is shown. The gain with the same dependence is displayed in subplot (f). We observe that the high gain in the early stage of the iteration is responsible for the rapid convergence of our method.}\label{Ex5}
	\end{center}
\end{figure}

The corresponding results are depicted in Fig.~\ref{Ex5}. As subplot~(a) reveals, contrary to the other cases, the obtained results for the inverse functions differ. For the initial value chosen, Newton's method breaks down for $y\geq3.32$. No such behavior is observed for the gMGF approach. A closer investigation of Newton's method for $y=3.33$ reveals that $|x_n\New|>|x_{n-1}\New|$ such that instead of converging to the numerical solution $\zeta\approx1.55529$ the iteration diverges.

From an inspection of the number of iterations (c.f. subplot~(b)), we observe that $I\gMGF<I\New$ for all considered values of $y$ within the interval $[0.1,3.32]$. Whereas the gMGF iteration exhibits a moderate behavior with $I\gMGF\leq11$, we observe particularly high values for $I\New$ around $y=1.88$ as well as $y=3.15$. For the first of these values, this can be traced back to a negative gain $g\New_1(y=1.88)\approx-0.947131$ of high magnitude during the first iteration. Interestingly, for $y=3.15$ we have $g_1\New(y=3.15)\approx0.300583$, whereas in the second iteration step we find a loss of $g_2\New(y=3.15)\approx-0.921243$. As observed during the previous examples, such a big loss regularly causes a high number of iterations to reach convergence.

Turning towards the number of function evaluations in subplot~(c), we observe, that for small $y$, $E\gMGF>E\New$ in spite of $I\gMGF<I\New$. This is due to the fact, that these points correspond to high values of $|\varkappa(x_n)|$. As the evaluation of $\mathfrak{H}_\varkappa$ alone requires $\varkappa$ function evaluations, the strong increase of $E\gMGF$ for small $y$ is explained. The situation is reversed however for $1.17\leq y\leq 2.4$, where, as discussed before, Newton's method exhibits a large increase in $I\New$. Consequently within this range we have $E\gMGF\leq E\New$. To the contrary, for $2.41\leq y\leq3$ we find $E\gMGF>E\New$. Beyond $y=3$ this situation is reversed again up to the value, where Newton's iteration eventually breaks down.

The number of function evaluations clearly influences the computing times, as can be seen from subplot (d). In the region where $E\gMGF>E\New$ holds, we find $t\gMGF>t\New$. Overall $t\gMGF$ remains quite moderate for the whole range of investigation which spans over five orders of magnitude. There is only a slight increase even towards its upper end. In comparison to this, Newton's method again reveals to be far more sensitive with respect to the distance $|\zeta-x_0|$.

The convergence behavior is elucidated in subplot~(e), where $y=1.5$ was chosen. We clearly observe, that the first iteration step within the gMGF approach takes approximately the time needed for three evaluations of Newton's method. However, we also clearly see, that, due to $|e\New_1|>|e\New_0|$, our method is superior in the investigated case for all iteration steps.

From subplot~(f) the loss $g\New_1\approx-0.47335$ is displayed. Subsequently Newton's method recovers very slowly, i.e. $|g\New_n|<|g\New_1|$ with $n\geq1$. The maximum gain is given for the tenth iteration step with $g\New_{10}\approx0.15514$. To the contrary, for the gMGF approach we find $g\gMGF_n\geq0$ and $g\gMGF_{n+1}<g\gMGF_n$, where $g\gMGF_1\approx0.62536$.

Tab.~\ref{tabular:T_MSL5} resembles some further data towards the convergence behavior of both algorithms, where $y=1.5$ was chosen. Again the negative value for $\varrho\New_2$ signals a loss connected to the initial iteration, which comparing $|x_1-x_0|$ to the corresponding value for higher $n$ displays a big magnitude. We see, that for $n\geq2$ the method converges to the actual solution. However the initial deterioration makes 17 iterations necessary to reach the convergence goal. Moreover, we observe that $\varrho\New_n$ strictly increases and again displays the theoretical value towards the end of the procedure. The initial iteration of the gMGF instead exhibits a gain which is signaled by both $\varrho_2\gMGF>0$ and the magnitude of the first figure in the gMGF section of the table. Overall a rapid convergence results such that after four iterations the set goal is reached. Again, we see that $\varrho_n\gMGF$ is less regular as the corresponding quantity within Newton's method. Finally, comparing the CPU times of both approaches, we see, that the gMGF approach outperforms Newton's method.
\begin{figure}[t]
	\begin{center}
		\includegraphics[width=.99\linewidth]{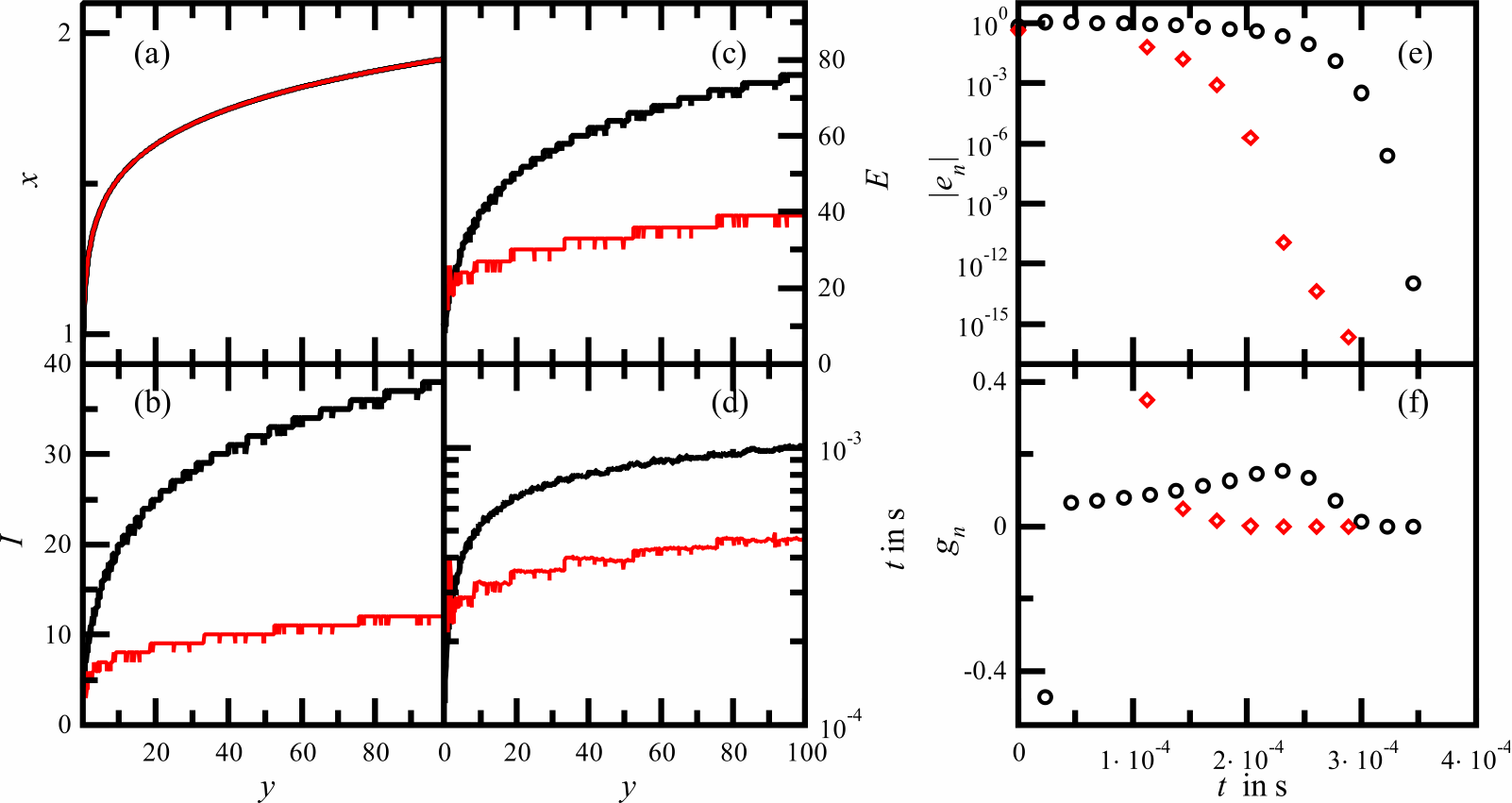}
		\caption{Comparison of results for \textbf{Example~\ref{Exam6}} as obtained from Newton's method (black curve) and our approach (red curve): subplot (a) displays the inverted function $x=x(y)$, whereas (b) shows the number of iterations needed for convergence. The number of function evaluations connected to this can be seen from subplot (c). In subplot (d) the computing time until convergence is shown. Especially when the starting value deviates strongly from the actual solution, our approach outperforms Newton's method. Subplots (e) and (f) show criteria to compare the convergence behavior of Newton's method and our approach for \textbf{Example~\ref{Exam6}} with $y=5$.  The circles display these quantities as obtained from Newton's method, whereas the diamonds are calculated from the approach suggested in this publication. In subplot (e) the residual as a function of the respective computing time is shown. The gain with the same dependence is displayed in subplot (f). We observe that the high gain in the early stage of the iteration is responsible for the rapid convergence of our method.}\label{Ex6}
	\end{center}
\end{figure}
\begin{table}[t]
	\centering
	\caption{Numerical results for~\textbf{Example 5:} $y=-\frac{1}{x}+\frac{1}{\sqrt{x}}+\frac{3}{20}x^{10}$ with $y=7$ and $x_0=1$}
	\label{tabular:T_MSL6}
	\begin{tabular}[!ht]{c c c c c c c c c c c c c c c}
	\toprule[2pt]
		Newton & $n=1$ & $n=2$ & $n=3$ & $n=4$ & $n=5$ & $n=6$&  $n=7$ &$n=8$\\
		\midrule
		$|x_n-x_{n-1}|$ & $3.4$ & $0.44$ & $0.40$ & $0.36$ & $0.32$ & $0.29$& $0.26$ & $0.23$\\
		$\varrho_n$ &  & $-0.088$ & $1.06$ & $1.08$ & $1.09$ & $1.11$& $1.14$ & $1.19$ \\\midrule
		& $n=9$ & $n=10$ & $n=11$ &$n=12$ &$n=13$ & $n=14$&$n=15$&$n=16$\\\midrule
		$|x_n-x_{n-1}|$ & $0.21$ & $0.18$ & $0.14$ & $0.091$ &$0.033$& $3.8\cdot10^{-3}$&$4.6\cdot10^{-5}$&$6.6\cdot10^{-9}$\\
		$\varrho_n$  & $1.25$ & $1.34$ & $1.47$ & $1.65$ & $1.84$ & $1.97$&$2.00$&$2.00$\\\midrule
		& &  & & & & &  $I$ & CPU time\\\midrule
		 &  &  &  &  & & & $16$&$4.31\cdot10^{-4}$ s\\\midrule\midrule
		gMGF & $n=1$ & $n=2$ & $n=3$ & $n=4$ & $n=5$ & $n=6$ & $I$ & CPU time\\
		\midrule
		$|x_n-x_{n-1}|$\ \  & $0.38$ & $0.12$ & $0.029$ & $2.8\cdot10^{-3}$ & $2.4\cdot 10^{-5}$ & $1.8\cdot 10^{-9}$ & $6$ & $2.68\cdot 10^{-4}$ s\\
		$\varrho_n$ &  & $0.590$ & $2.39$ & $1.96$ & $2.00$ & $2.00$ & $2.00$\\
		\bottomrule[2pt]
	\end{tabular} 
\end{table}
\begin{table}[t]
	\centering
	\caption{Numerical results for~\textbf{Example 6:} $y=x^9+x^7+x^2$ with $y=5$ and $x_0=0.2$}
	\label{tabular:T_MSL7}
	\begin{tabular}[!ht]{c c c c c c c c c c c c c c c}
	\toprule[2pt]
		Newton & $n=1$ & $n=2$ & $n=3$ & $n=4$ & $n=5$ & $n=6$&  $n=7$ &$n=8$\\
		\midrule
		$|x_n-x_{n-1}|$ & $12$ & $1.4$ & $1.2$ & $1.1$ & $0.98$ & $0.88$& $0.78$ & $0.69$\\
		$\varrho_n$ &  & $-0.0505$ & $1.01$ & $1.02$ & $1.02$ & $1.02$& $1.02$ & $1.03$ \\\midrule
		& $n=9$ & $n=10$ & $n=11$ &$n=12$ &$n=13$ & $n=14$&$n=15$&$n=16$\\\midrule
		$|x_n-x_{n-1}|$ & $0.62$ & $0.55$ & $0.49$ & $0.43$ &$0.39$& $0.34$&$0.31$&$0.28$\\
		$\varrho_n$  & $1.03$ & $1.04$ & $1.04$ & $1.05$ & $1.06$ & $1.07$&$1.09$&$1.12$\\\midrule
		& $n=17$ & $n=18$ & $n=19$ &$n=20$ &$n=21$ & $n=22$&$n=23$&$n=24$\\\midrule
		$|x_n-x_{n-1}|$ & $0.24$ & $0.22$ & $0.19$ & $0.16$ &$0.12$& $0.069$&$0.019$&$1.2\cdot10^{-3}$\\
		$\varrho_n$  & $1.15$ & $1.20$ & $1.27$ & $1.37$ & $1.52$ & $1.72$&$1.90$&$1.98$\\\midrule
		& $n=25$ & $n=26$  &  & & & &  $I$ & CPU time\\\midrule
		$|x_n-x_{n-1}|$ &  $4.7\cdot10^{-6}$ & $7.1\cdot10^{-11}$ & & & & & $26$&$3.63\cdot10^{-4}$ s\\\midrule
		$\varrho_n$& $1.99$ & $2.00$ &  &  & & &\\\midrule\midrule
		gMGF & $n=1$ & $n=2$ & $n=3$ & $n=4$ & $n=5$ & $n=6$ & $n$=7 \\
		\midrule
		$|x_n-x_{n-1}|$\ \  & $0.40$ & $0.41$ & $0.094$ & $0.019$ & $1.2\cdot 10^{-3}$ & $4.7\cdot 10^{-6}$ & $6.8\cdot10^{-11}$ & \\
		$\varrho_n$ &  & $3.07$ & $0.687$ & $2.18$ & $1.98$ & $2.00$ & $2.00$\\\midrule
		& & & & & & & $I$ & CPU time\\\midrule
		& & & & & & & $7$ & $2.74\cdot10^{-4}$ s\\
		\bottomrule[2pt]
	\end{tabular} 
\end{table}
\begin{ex}\label{Exam6}
	The next function, which we like to invert is given by
	\begin{equation}
		y=-\frac{1}{x}+\frac{1}{\sqrt{x}}+\frac{3}{20}x^{10},
	\end{equation}
	where we investigate the range $y\in[0.1,100]$ with the step size $\delta y=0.1$. As the initial value, we take $x_\mathrm{0}=1$.
\end{ex}
In Fig.~\ref{Ex6} results connected to \textbf{Example~\ref{Exam6}} are displayed. As can be seen from subplot (a) in the present case both methods succeed in inverting the function $y$ in the considered range. From subplot~(b) we observe that in order to reach this result, Newton's method requires more iterations then the iterative gMGF approach in the whole range taken into account. This is slightly different as far as the number of function evaluations are concerned. In fact for $y\in[0.1, 1.8]$ we have $E\New\leq E\gMGF$. Above this value the gMGF algorithm performs better than Newton's method in this respect. Turning towards the computing time (c.f. subplot~(d)), the same tendency is observed. However, due to details in the evaluation of $\mathfrak{H}_\varkappa$, we find that the interval, where Newton's method outperforms the gMGF approach is extended to $[0.1, 1.9]$. Within the remainder of the considered range, the gMGF method clearly offers advantages as far as computing time is concerned. At the upper end of the investigated range, the computation time of Newton's method is cut in half by the iterative gMGF approach. We note, that $\Delta t(y)\sim\ln y$ persists beyond the upper boundary of the shown interval.

A more detailed analysis of the convergence behavior by means of $|e^{(\gamma)}_n|$ displayed in subplot~(e) for $y=5$ exhibits the same tendency which was reported for the preceding examples: The evaluation of a single iteration within the gMGF requires a longer computing time than Newton's method. However, as $|e_1\New(5)|>|e_0\New(5)|$, Newton's method suffers from the sensitivity with respect to the initial value and in total it takes seven iteration steps to compensate the error resulting from the initial calculation, i.e. $|e_l\New(5)|>|e_0\New(5)|$ for $l\in\{1,2,\dots,7\}$. Overall, 15 iterations are required for convergence within Newton's method. No such behavior is observed within the gMGF approach and the comparatively long computing times per iteration are rewarded with a smaller number of iteration steps ($I\gMGF(5)=7$), thus resulting in an altogether smaller $t\gMGF(5)\approx0.28852$ ms to meet the demanded precision.

Finally, the gain of the two methods is depicted in subplot~(f). Here, the loss of Newton's method within the first iteration step is clearly displayed. Concerning its magnitude, it exceeds the gain obtained from the second step by over a factor of five. Contrary to this, the gMGF approach displays a positive gain of $g_1\gMGF(5)\approx0.35094$, followed by the aforementioned rapid convergence.

For $y=7$ some additional insights into the convergence behavior are collected in Tab.~\ref{tabular:T_MSL6}. As for the previous examples, the number of iterations within Newton's method exceeds the one connected to the gMGF algorithm. In total, we find $I\New=16$. Again we observe an initial loss displayed in the COC for $n=2$ connected to a comparatively large magnitude $|x_1-x_0|$. The consecutive values in this line turn out to be much smaller and a slow convergence is observed. The COC strictly increases with $n$ reaching the theoretical value towards the end of the calculation. Altogether, the gMGF requires six iterations for the example under study. The quantity $|x_n-x_{n-1}|$ turns out be strictly decreasing, whereas $\varrho\gMGF$ changes significantly for the first iterations until for $n\geq4$ the COC also converges to 2.00.
\begin{figure}[t]
	\begin{center}
	\includegraphics[width=.99\linewidth]{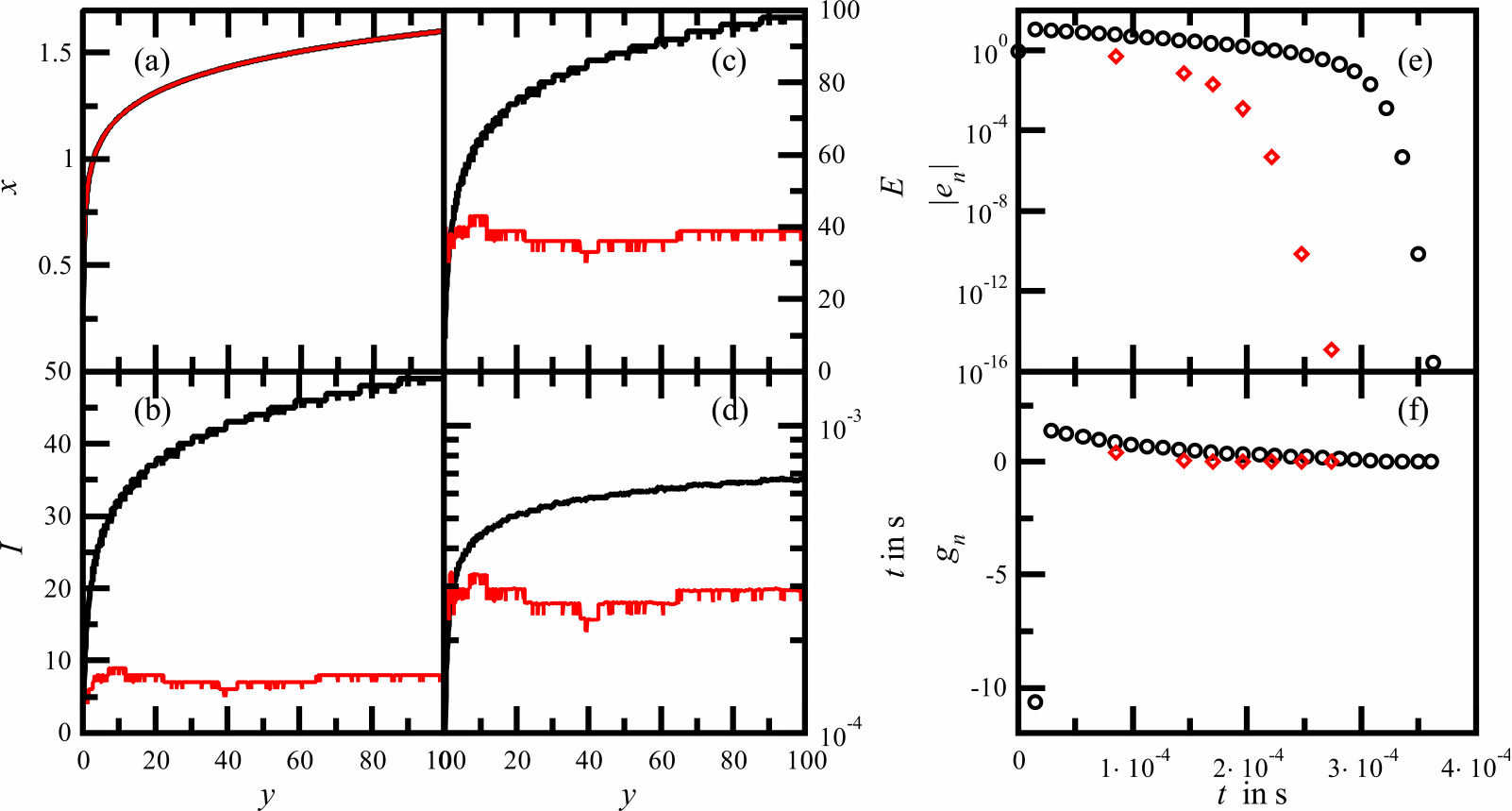}
	\caption{Comparison of results for \textbf{Example~\ref{Exam8}} as obtained from Newton's method (black curve) and our approach (red curve): subplot (a) displays the inverted function $x=x(y)$, whereas (b) shows the number of iterations needed for convergence. The number of function evaluations connected to this can be seen from subplot (c). In subplot (d) the computing time until convergence is shown. Especially when the starting value deviates strongly from the actual solution, our approach outperforms Newton's method. Subplots (e) and (f) show criteria to compare the convergence behavior of Newton's method and our approach for \textbf{Example~\ref{Exam8}} with $y=5$. The circles display these quantities as obtained from Newton's method, whereas the diamonds are calculated from the approach suggested in this publication. In subplot (e) the residual as a function of the respective computing time is shown. The gain with the same dependence is displayed in subplot (f).}\label{Ex8}
	\end{center}
\end{figure}
\begin{ex}\label{Exam8}
Next we choose the polynomial
	\begin{equation}
		y=x^9+x^7+x^2,
	\end{equation}
which we want to invert for $y\in[0.1,99.1]$ with $\delta y=0.1$ and the initial value $x_0=0.2$.
\end{ex}
Results from the iterations are depicted in Fig~\ref{Ex8}. Again we observe from subplot~(a), that both methods converge to $x(y)$ in the considered range for $y$. For the number of iterations, we have $I\gMGF(y)<I\New(y)$ for all values taken into consideration as can be seen from subplot~(b). In particular, we find a strong increase of $I\New(y)$ with increasing $y$ for $y\leq10$. The maximum number of iterations within Newton's method amounts to 49. This is to be compared with he corresponding quantity from the gMGF algorithm which is given by 9 occurring around $y=10$.

A somewhat different behavior is observed for the number of function evaluations displayed in subplot~(c). Here we find $E\gMGF(y)>E\New(y)$ for $y\leq2.1$ whereas the relation between the function evaluations of both methods is reversed within the rest of the considered range.

As a result of this also the computing time, displayed in subplot~(d), reveals the same tendency. While $t\New(y)$ increases with $y$, the corresponding quantity within the gMGF approach shows only little variation. Analogously to the number of function evaluations, within the given resolution $\delta y$, we find $t\gMGF(y)<t\New(y)$ for $y\geq2.2$. For $y=100$ the computation time of the iterative gMGF approach amounts to half the one needed for convergence within Newton's method. As in the previous example, $\Delta t(y)$ exhibits a logarithmic tendency which also occurs for $y>100$.

From subplot~(e) we again observe that for $y=5$ the magnitude of the residual in the first iteration step increases within Newton's method. Such a behavior is not observed for the gMGF approach where we find $|e_{n+1}\gMGF|<|e_{n}\gMGF|$. Also, we find the same tendency as seem from our previous examples, i.e., that the single iteration steps require more time within the gMGF approach as compared to Newton's method. In the present case, we find, that the first iteration within the gMGF approach takes approximately the time of six steps within Newton's method. However, we observe a more rapid convergence for the gMGF approach due to the reduced number of required iterations.

As far as the gain, displayed in subplot~(f), is concerned, as in most previous results, we find an initial loss for Newton's method which is absent for any of the steps of our algorithm. In the present case, the magnitude of $g\New_{1}$ is particularly high. For higher values of $y$ and thus $|\zeta-x_0|$, the initial loss exhibits an even higher magnitude, such that the benefits of the gMGF method increases accordingly, as can be seen from subplot~(d).

From Tab.~\ref{tabular:T_MSL7} some details about the convergence behavior can be obtained also for $y=5$. Starting from an initial loss, detectable from $\varrho_2\New$, Newton's method takes 26 iteration steps until convergence is reached. We see, that for $n\geq2$ $|x_n-x_{n-1}|$ slowly decreases, while $\varrho\New_n$ shows an increase which eventually tends towards the theoretical value 2. Within the seven iterations required for convergence within the gMGF approach, no loss is detected. The COC varies non-monotonically for the first iterations until for $n\geq5$ a tendency towards the theoretical value is observed.
\begin{ex}\label{Exa7} Our next example aims at inverting the function
	\begin{equation}\label{BBR}
		y=e^{-x}+\frac{x}{5},
	\end{equation}
	with $y\in[1,70]$ and the step size $\delta y=0.1$. As the initial value, we choose $x_0=1$.
\end{ex}
\begin{figure}[t]
	\begin{center}
	\includegraphics[width=.99\linewidth]{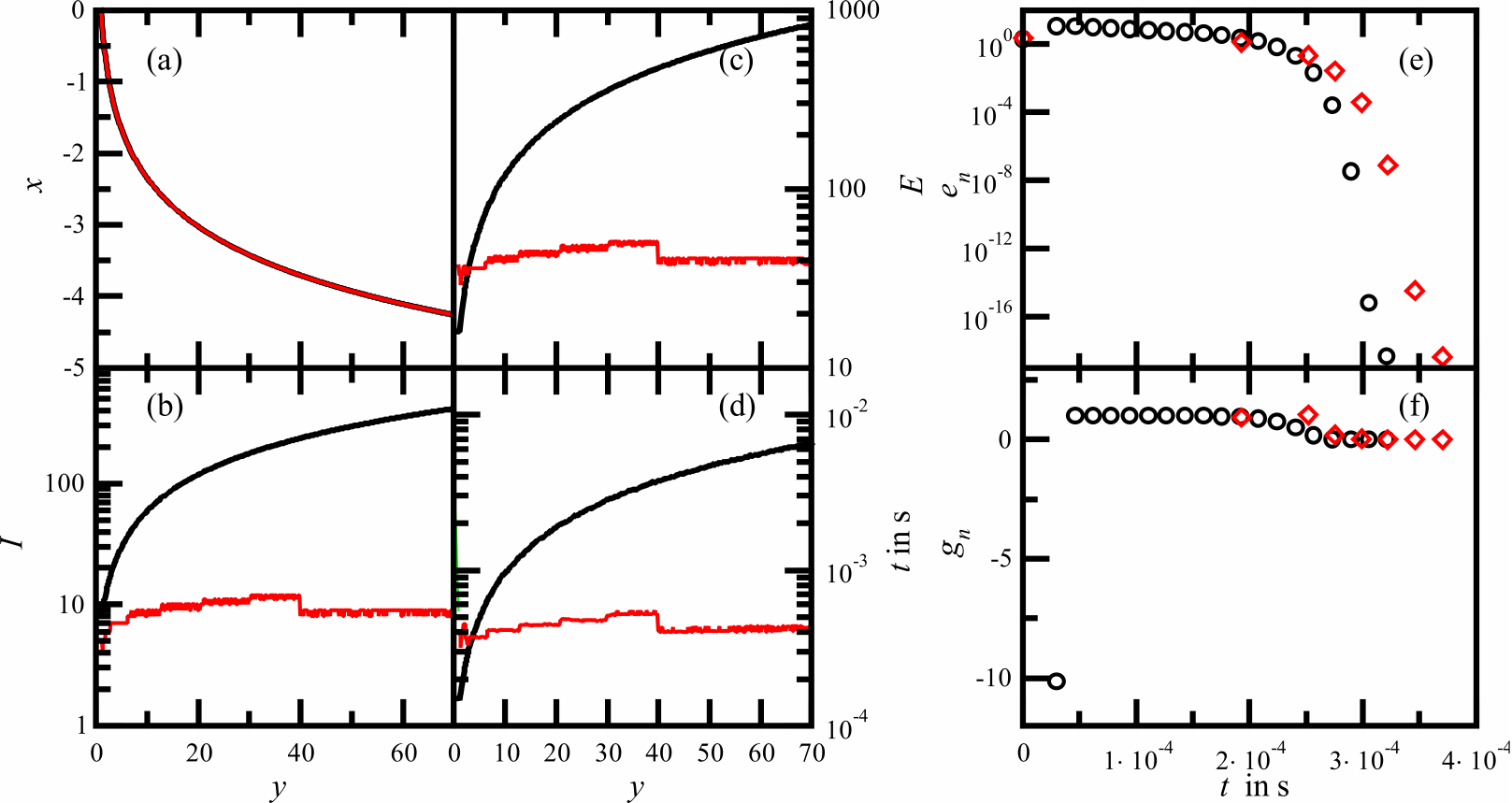}
	\caption{Comparison of results for \textbf{Example~\ref{Exa7}} as obtained from Newton's method (black curve) and our approach (red curve): subplot (a) displays the inverted function $x=x(y)$, whereas (b) shows the number of iterations needed for convergence. The number of function evaluations connected to this can be seen from subplot (c). In subplot (d) the computing time until convergence is shown. Especially when the starting value deviates strongly from the actual solution, our approach outperforms Newton's method. Subplots~(e) and~(f) show criteria to compare the convergence behavior of Newton's method and our approach for \textbf{Example~\ref{Exa7}} with $y=3$. The circles display these quantities as obtained from Newton's method, whereas the diamonds are calculated from the approach suggested in this publication. In subplot (e) the magnitude of the residual as a function of the iteration step is shown. The gain is displayed in subplot (f).}\label{Ex7}
	\end{center}
\end{figure}
\begin{table}[t]
	\centering
	\caption{Numerical results for~\textbf{Example 7:} $y=e^{-x}+\frac{x}{5}$ with $y=3$ and $x_0=1$.}
	\label{tabular:T_E7}
	\begin{tabular}[!ht,width=.98\linewidth]{c c c c c c c c c c c c c c c}
	\toprule[2pt]
		Newton & $n=1$ & $n=2$ & $n=3$ & $n=4$ & $n=5$ & $n=6$ & $n=7$ \\
		\midrule
		$|x_n-x_{n-1}|$ & $14$ & $1.0$ & $1.0$ & $1.0$ & $1.0$ & $1.0$ & $ 1.0$ \\
		$\varrho_n$ &  & $-4.88\cdot10^{-2}$ & $1.09$ & $1.10$ & $1.11$ & $1.13$ \\\midrule
		&$n=8$ &  $n=9$ & $n=10$ & $n=11$ & $n=12$ & $n=13$ & $n=14$ \\\midrule
		$|x_n-x_{n-1}|$ & $1.0$ & $0.99$ & $0.98$ & $0.96$ & $0.90$ & $0.76$ & $ 0.51$  \\
		$\varrho_n$ & $1.14$ & $1.17$ & $1.20$ & $1.25$ & $1.31$ & $1.40$ \\\midrule
		& $n=15$ & $n=16$ &  $n=17$ & $n=18$ & $n=19$ &   $I$ & CPU time \\\midrule
		$|x_n-x_{n-1}|$ & $0.19$ & $0.022$&  $2.6\cdot10^{-4}$ & $3.6\cdot10^{-8}$ & $6.7\cdot10^{-16}$   &    19 & $3.21\cdot10^{-4}$ \\
		$\varrho_n$ & $1.52$ & $1.68$ & 1.85 & $1.97$ & $2.00$  &  &  &  \\\midrule\midrule
		gMGF & $n=1$ & $n=2$ & $n=3$ & $n=4$ & $n=5$  & $n=6$ & $n=7$\\
		\midrule
		$|x_n-x_{n-1}|$\ \  & $0.92$ & $1.0$ & $0.24$ & $2.7\cdot10^{-2}$ & $3.9\cdot 10^{-4}$ & $7.9\cdot10^{-8}$ & $3.3\cdot10^{-15}$\\
		$\varrho_n$ &  & $3.20$ & $1.18$ & $2.05$ & $2.00$ & $2.00$ & $2.00$\\\midrule
		 &  &  &  &  &   & $I$ & CPU time\\
		\midrule
		 & & & & & & $7$ & $3.70\cdot 10^{-4}$ s\\
		\bottomrule[2pt]
	\end{tabular} 
\end{table}
Fig.~\ref{Ex7} reveals the results for the example in question. As before, subplot~(a) demonstrates that both methods reproduce the inverse function $x(y)$ accurately in the whole range under consideration. The number of iterations (c.f. subplot~(b)) reveal, that $I\gMGF(y)<I\New(y)$ holds. Due to the calculation of $\varkappa$ and $\mathfrak{H}_\varkappa$, the number of function evaluations for Newton's method, however, is smaller than within the gMGF framework for $y<2.8$. For $y\geq2.9$, we find $E\gMGF(y)\leq E\New(y)$ where equality holds for $y=2.9$.

Turning towards the computing time displayed in subplot~(d), we observe, that these graphs follow the same trend as the number of function evaluations. In particular, we find, that for $y<3.7$ the relation $t\New(y)<t\gMGF(y)$ holds whereas the computing time within the gMGF approach turns out to be smaller than by means of Newton's method for $y\geq3.8$. While $t\gMGF\leq0.65695$ ms remains relatively moderate and does not change strongly within the considered range, the maximum value within Newton's iteration reads $t\New=4.4039$ ms. The computing time increases linearly with increasing $y$, whereas no such tendency is observed for the gMGF approach. At the upper boundary of the interval, the computation time within Newton's method is twenty times as high as for the iterative gMGF approach.

To investigate the convergence behavior for this example, we study the residuals and the gains for a value of $y=3$, where, as far as computation time is concerned, Newton's method outperforms our algorithm. Still, we observe, that within the first step of Newton's method, the residual increases as for most of the previous examples. No such behavior is seen for the gMGF algorithm. Instead, $|e_n\gMGF|$ decreases with increasing $n$. Moreover, we observe, that the single iteration steps within our method take considerably longer than within Newton's algorithm. On the other hand, the decrease in magnitude of the residuals turn out to be larger within the gMGF scheme. Interestingly in the first half of the graph, we see that our iteration outperforms Newton's method. Therefore, also for this example we find evidence, that the gMGF approach is particularly suitable for applications with a comparatively low precision goal, for in this case the main advantage, of the method, namely the low number of necessary iterations comes into full play.

Turning towards the gain in subplot~(f), as before, we recognize a big loss within Newton's method for the first step which is compensated during the subsequent iterations. Interestingly, also within the gMGF method for the example, the first step does not exhibit maximum gain. Instead, we find that this is reached at the second iteration step with $g\gMGF_2\approx1.039$.

In Tab.~\ref{tabular:T_E7} we compare the quantities $|x_n-x_{n-1}|$ as well as the COC as obtained from the two methods during the iteration for $y=3$. Newton's method again exhibits a negative value for the COC for $n=2$. The initial loss connected to it causes the algorithm to make only slow progress, as can be seen from the almost constant tendency of $|x_n-x_{n-1}|$ between $n=2$ and $n=10$. We would like to point out, that $\varrho_n\New$ is also only slightly increasing for these iteration steps. Thereupon, the considered quantities vary more rapidly and the COC tends towards the corresponding theoretical value until convergence is reached with a total of 19 iteration steps. In contrast to the discussed behavior, again the gMGF approach shows rapid convergence from the first iteration step and a total of 7 iteration steps are sufficient to reach convergence. Albeit this difference, we observe, that the CPU time of Newton's method outperforms the gMGF algorithm for the value under investigation.
\section{Applications}\label{Apps}
In this section, we want to compare the performance of our approach to Newton's method for real world applications. In doing so, we have intentionally chosen objectives, where the initial value for the iteration is hard to determine, while the process under study requires an evaluation within a dependable computing time.
\subsection{Calculation of the principal branch of the Lambert $W$ function}\label{Exam3}
The principal branch of the Lambert $W$ function results from resolving the equation
\begin{equation}\label{LambertW}
	y=x e^{x}
\end{equation}
for $x$ with $y>-\exp\{-1\}$.
\begin{figure}[t]
	\begin{center}
	\includegraphics[width=.99\linewidth]{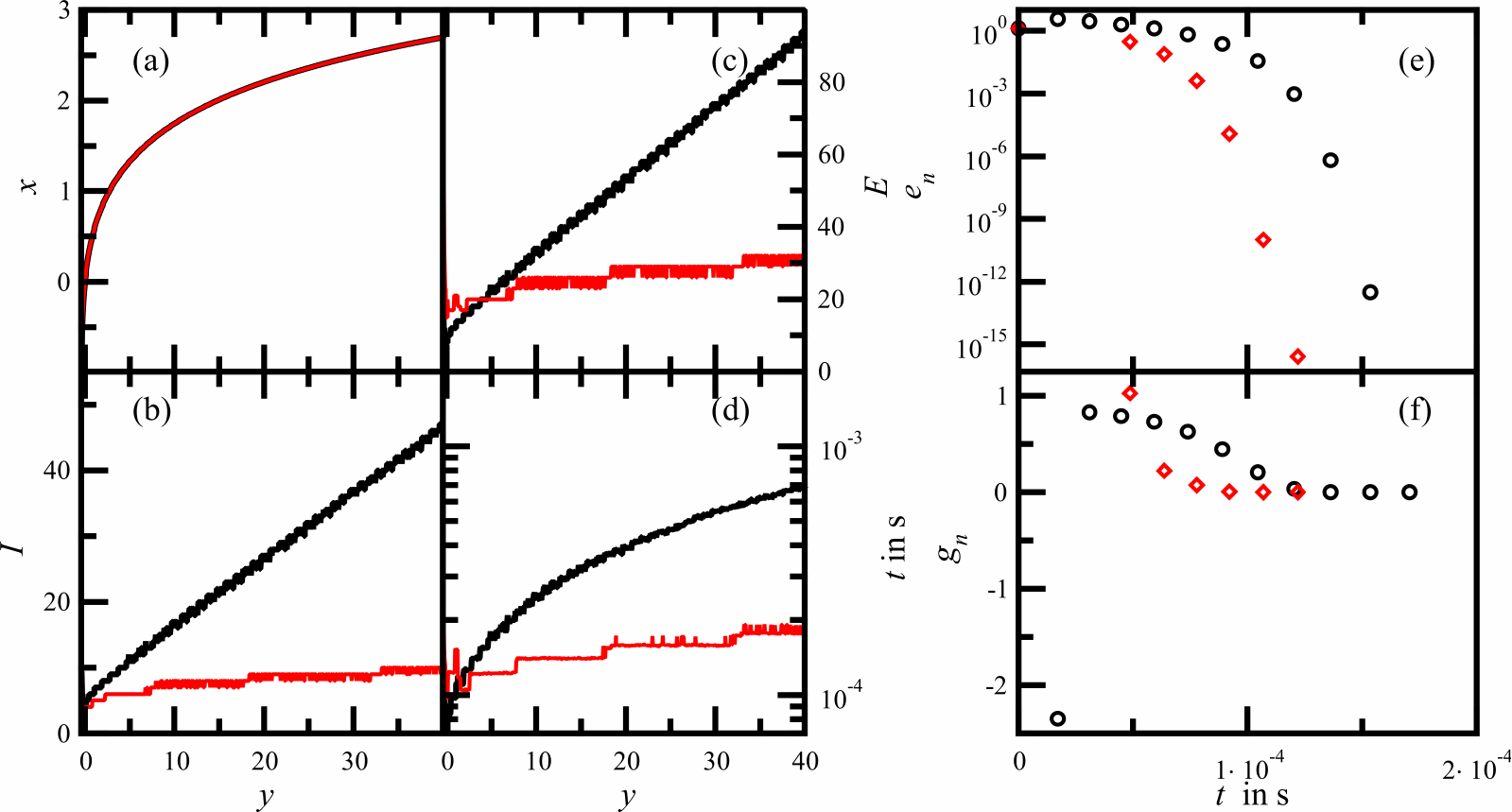}
	\caption{Comparison of results for solving Eq.~(\ref{LambertW}) as obtained from Newton's method (black curve) and our approach (red curve): subplot (a) displays the inverted function $x=x(y)$, whereas (b) shows the number of iterations needed for convergence. The number of function evaluations connected to this can be seen from subplot (c). In subplot (d) the computing time until convergence is shown. Especially when the starting value deviates strongly from the actual solution, our approach outperforms Newton's method. Criteria to compare the convergence behavior of Newton's method and our approach for  solving Eq.~(\ref{LambertW}) with $y=5$ are shown in subplots~(e) and~(f). The circles display results as obtained from Newton's method, whereas the diamonds are calculated from the approach suggested in this publication. In subplot (e) the magnitude of the residuals as a function of the computing time is shown. The gain is displayed in subplot (f).}\label{Ex3}
	\end{center}
\end{figure}
This function covers a broad spectrum of applications including solutions for delay differential equations~\cite{Wright}, control algorithms for systems with memory~\cite{Sun2008}, as well as the current-voltage characteristics of solar cells~\cite{JAIN2005197, KAPOOR2006120, JAIN2004269} and Schottky diodes~\cite{JUNG200957}. It also is used in the field of fluid dynamics, as the Fanning friction factor for turbulent conditions and also the impact of the exhaust turbo charger on the formation of nanoscale particulate matter for a combustion engine can be expressed by the Lambert W function ~\cite{MORE20065515, CUCCHI2016429}. For these reasons an efficient implementation might be of particular interest for process engineering and automotive applications. Subsequently, we shall see, that solving Eq.~(\ref{LambertW}) can significantly be sped up by using our approach instead of Newton's method.

As was the case in the previous section, rather than to consider different initial values we choose to vary the function value as follows $y\in[-0.367,39.983]$. In doing so, we use the step size $\delta y=0.05$. For the initial value we use $x_0=0$.

In Fig.~\ref{Ex3} the corresponding results to this objective are displayed. The Lambert $W$ function is clearly reproduced by both methods taken into consideration as can be seen from subplot~(a). A closer investigation of the number of iterations in subplot~(b) reveals, that both methods are at par at $y=-0.367$ with $I\New(-0.367)=I\gMGF(-0.367)=3$, whereas for the rest of the values considered, we have $I\New(y)>I\gMGF(y)$. Also, we recognize an essentially linear increase of $I\New$ with $y$, whereas $I\gMGF$ reveals a flat plateau-like increase. 

The same tendency is observed for the number of function evaluations (c.f. subplot~(c)), where, with the same reasons as mentioned in the previous examples $E\New(y)\leq E\gMGF(y)$ for $y\in[-0.367,4.583]$. Within the rest of the considered range, we find $E\gMGF(y)<E\New(y)$. Interestingly, considering the computing time, the interval where Newton's method outperforms the gMGF algorithm is reduced to $y\in[-0.367,1.533]$. This is due to the fact, that the additional operations processed at these points within the iterative gMGF are connected to a low computational cost. Whereas the computing time for Newton's method increases essentially linearly with increasing $|y-f(x_0)|$, the moderate plateau-shaped increase already observed for the number of iterations is also displayed here for the gMGF approach.

The investigations undertaken for the convergence behavior in this case refer to $y=5$. From subplot~(e), we observe, as in the examples discussed above, that the single iterations within the gMGF algorithm require more time than those of Newton's method. However, the magnitudes of the residuals reveal, that our method performs more efficient during the whole computation for the considered value. In particular, we again observe, that within the first iteration step, $|e\New_1|>|e\New_0|$, a behavior which is absent within the gMGF approach. From subplot~(f), we see, that the initial loss caused for Newton's method necessitates a total number of 11 iterations to meet the precision goal. To the contrary, due to the uniform convergence of the gMGF approach, this aim is reached after 6 steps.

In Tab.~\ref{tabular:T_MSL3} some further details on the convergence behavior of the two methods for $y=5$ are displayed. From $\varrho_2\New<0$ we may again observe the loss which results from the first iteration step within Newton's method. Thereupon the method converges with an overall CPU time of $1.71\cdot10^{-4}$ s. The COC exhibits a strictly increasing tendency and it reaches the theoretical value towards the end of the overall procedure. The corresponding quantity for the gMGF method varies stronger over the six iterations. However, we also observe, that the theoretical value is reached. In particular, we find, that again the rapid convergence in the early iteration steps makes this algorithm overall less time-consuming than Newton's method, albeit the single iterations within the gMGF approach require more time.
\begin{table}[t]
	\centering
	\caption{Numerical results for solving Eq.~(\ref{LambertW}) for $x$ with $y=5$ and $x_0=0$}
	\label{tabular:T_MSL3}
	\begin{tabular}[!ht]{c c c c c c c c c c c c c c c}
	\toprule[2pt]
		Newton & $n=1$ & $n=2$ & $n=3$ & $n=4$ & $n=5$ & $n=6$  \\
		\midrule
		$|x_n-x_{n-1}|$ & $5.0$ & $0.83$ & $0.79$ & $0.73$ & $0.63$ & $0.45$ \\
		$\varrho_n$ &  & $-0.248$ & $1.28$ & $1.35$ & $1.46$ & $1.61$ \\\midrule
		&  $n=7$ & $n=8$ & $n=9$ & $n=10$ & $n=11$ & & $I$ & CPU time & &\\\midrule
		$|x_n-x_{n-1}|$ & $0.21$ & $0.036$ & $9.7\cdot10^{-4}$ & $6.7\cdot10^{-7}$ & $3.2\cdot10^{-13}$ & &$11$&$1.71\cdot10^{-4}$ s\\
		$\varrho_n$ & $1.79$ & $1.94$ & $1.99$ & $2.00$ & $2.00$ \\\midrule\midrule
		gMGF & $n=1$ & $n=2$ & $n=3$ & $n=4$ & $n=5$ & $n=6$ & $I$ & CPU time\\
		\midrule
		$|x_n-x_{n-1}|$\ \  & $1.03$ & $0.38$ & $7.3\cdot10^{-2}$ & $4.1\cdot10^{-3}$ & $1.2\cdot 10^{-5}$ & $1.0\cdot10^{-10}$ & $6$ & $1.22\cdot 10^{-4}$ s\\
		$\varrho_n$ &  & $0.916$ & $2.16$ & $1.99$ & $2.00$ & $2.00$\\
		\bottomrule[2pt]
	\end{tabular} 
\end{table}

\subsection{Flow-rate control in a bio reactor}
A large number of various products like food, chemicals, pharmaceuticals and agriculture stem from industrial bio processing. From this perspective the monitoring and controlling of cell growth alongside the bio catalytic and biological processes are of importance. For a stable and efficient process, the population of the biomass in a bio reactor is of interest. For instance in a chemostat, a certain flow rate $\dot V$ of nutrients required to synthesize the cells is transferred through a volume $V$. Particular high cell densities can be achieved by recycling the cells in the outlet stream back into the vessel. If we denote the recycle ratio by $\Theta$ and the factor by which the cells are concentrated before being fed back into the volume by $\Phi$, the cell concentration $c$ in such an open system can be modeled via
\begin{equation}\label{chemostat}
	\frac{dc}{dt}=\frac{\dot V}{V}\left[c_0-c\left(1+\Theta-\Phi\right)\right]+\left(\Gamma-\Delta\right)c.
\end{equation}
Here $\Gamma$ and $\Delta$ respectively resemble the specific growth and death rate of the cells, while $c_0$ resembles the initial concentration of the biomass~\cite{NAJAFPOUR200781}. Introducing the abbreviations $y=c/c_0$, $\lambda=\Gamma-\Delta$, $\nu=\Phi-\Theta$, and $x=\dot V/V$, which can be regarded as the volume specific flow rate, we may solve Eq.~(\ref{chemostat}), to read
\begin{equation}\label{Grow}
	y(t)=e^{[\lambda+(\nu-1)x]t}+\frac{x}{\lambda+(\nu-1)x}\left(e^{[\lambda+(\nu-1)x]t}-1\right).
\end{equation}
In order to reach the desired cell growth subject to the technical limitations of the process, it is necessary to maintain a certain value for the flow rate. The objective in this case is to calculate $x$ from Eq.~(\ref{Grow}). For the following simulation, we set $\lambda=0.8\ \textnormal{h}^{-1}$, $\nu=0$ and $t=10$ h. We investigate $y\in\left[1.3,8\right]$ with a step size of $\delta y=0.01$ and the initial value is chosen as $x_0=2$.
\begin{figure}[t]
	\begin{center}
	\includegraphics[width=.99\linewidth]{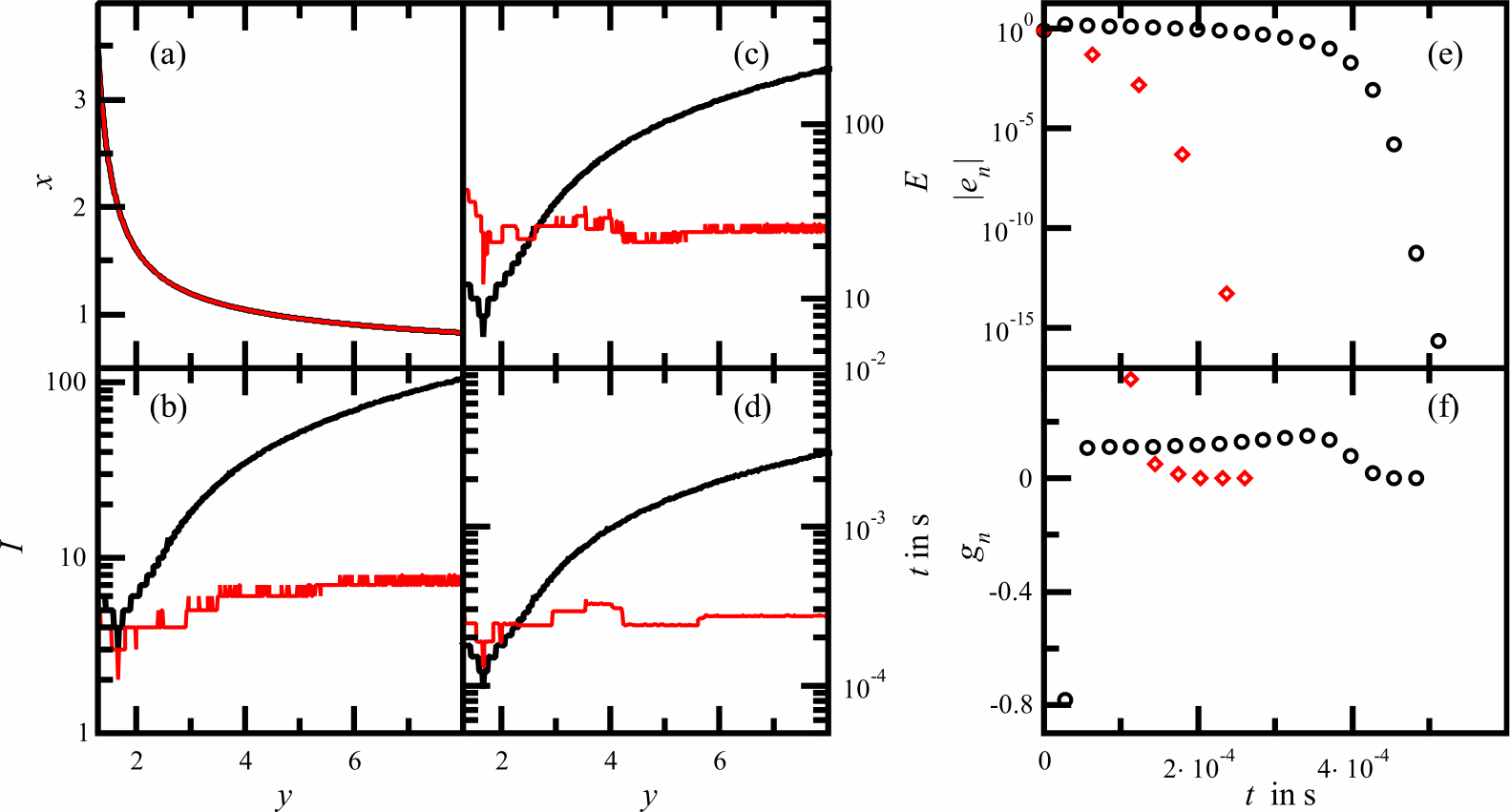}
	\caption{Comparison of results for solving Eq.~(\ref{Grow}) as obtained from Newton's method (black curve) and our approach (red curve): subplot (a) displays the inverted function $x=x(y)$, whereas (b) shows the number of iterations needed for convergence. The number of function evaluations connected to this can be seen from subplot (c). In subplot (d) the computing time until convergence is shown. Especially when the starting value deviates strongly from the actual solution, our approach outperforms Newton's method. Criteria to compare the convergence behavior of Newton's method and our approach for  solving Eq.~(\ref{Grow}) with $y=3$ are shown in subplots~(e) and~(f). The circles display results as obtained from Newton's method, whereas the diamonds are calculated from the approach suggested in this publication. In subplot (e) the magnitude of the residuals as a function of the computing time is shown. The gain is displayed in subplot (f).}\label{Gro}
	\end{center}
\end{figure}
\begin{table}[t]
	\centering
	\caption{Numerical results for solving Eq.~(\ref{Grow}) for $x$ with $\lambda=0.8$, $\nu=0$, $t=10$, $y=3$ and $x_0=2$}
	\label{tabular:Grow}
	\begin{tabular}[!ht]{c c c c c c c c c c c c c c c}
	\toprule[2pt]
		Newton & $n=1$ & $n=2$ & $n=3$ & $n=4$ & $n=5$ & $n=6$ & $n=7$ & $n=8$   \\
		\midrule
		$|x_n-x_{n-1}|$ & $2.4$ & $0.11$ & $0.11$ & $0.11$ & $0.12$ & $0.12$ & $0.12$ & $0.13$ \\
		$\varrho_n$ &  & $-105$ & $1.09$ & $1.10$ & $1.11$ & $1.13$ & $1.15$ & $1.18$\\\midrule
		& $n=9$ & $n=10$ & $n=11$ & $n=12$ & $n=13$ & $n=14$ & $n=15$ & $n=16$\\\midrule
		$|x_n-x_{n-1}|$  & $0.14$ & $0.14$ & $0.15$ & $0.14$ &$0.081$&$0.020$&$8.9\cdot10^{-4}$&$1.7\cdot10^{-6}$\\
		$\varrho_n$ & $1.23$ & $1.30$ & $1.40$ & $1.55$ & $1.72$ & $1.88$ & $1.98$ & $2.00$\\\midrule
		& $n=17$ & & & & & & $I$ & CPU time & &\\\midrule
		$|x_n-x_{n-1}|$ & $6.0\cdot10^{-12}$ & & &  &  & &$17$&$5.11\cdot10^{-4}$ s\\
		$\varrho_n$ & $2.00$ &  &  \\\midrule\midrule
		gMGF & $n=1$ & $n=2$ & $n=3$ & $n=4$ & $n=5$ &  & $I$ & CPU time\\
		\midrule
		$|x_n-x_{n-1}|$\ \  & $0.86$ & $0.055$ & $1.5\cdot10^{-3}$ & $5.0\cdot10^{-7}$ & $5.7\cdot 10^{-14}$ &  & $5$ & $2.47\cdot 10^{-4}$ s\\
		$\varrho_n$ &  & $1.31$ & $2.25$ & $1.99$ & $2.00$ & \\
		\bottomrule[2pt]
	\end{tabular} 
\end{table}

Results from this investigation are displayed in Fig.~\ref{Gro}. From subplot~(a), we observe that both, the result from Newton's method (black curve), as well as the gMGF approach capture the flow rate volume specific flow rate $x$ to high accuracy. Turning towards the number of iterations in subplot~(b), we find, that in the whole range $I\gMGF(y)<I\New(y)$ applies. While this quantities exhibits a linear tendency with increasing $y$ for Newton's method, for the gMGF algorithm, we again observe a slowly increasing plateau-like curve. While within the investigated interval $I\New$ grows up to 140, the number of iterations within the gMGF approach does not exceed eight.

Turning towards the number of functions evaluations, we find that $E\New(y)\leq E\gMGF(y)$ for $y\leq 2.7$. Above this value the gMGF algorithm outperforms Newton's method in this respect. In particular, we observe that $E\gMGF(y)$ stays essentially constant, other than the linear increase occurring within Newton's method.

The tendencies of the number of function evaluations are also resembled in the computation time until convergence, even though the interval, where the gMGF approach is inferior is reduced to $[1.3,2.31]$. In particular, we find that the linear increase with increasing $y$ for Newton's method is also seen for $t(y)$, whereas the computation time within our algorithm remains almost constant and fulfills $1.3074\cdot10^{-4}\ \textnormal{s}\leq t\gMGF\leq3.4748\cdot10^{-4}\ \textnormal{s}$ in the considered range. This has to be compared with the range that occurs within Newton's method, reading $9.4796\cdot10^{-5}\ \textnormal{s}\leq t\leq 2.9033\cdot10^{-3}\ \textnormal{s}$.

Subplot~(e) displays the magnitude of the residuals as a function of computation time for $y=3$. It reveals the same tendencies as reported in most previous examples: In the first step Newton's method causes a deterioration with respect to the initial vale. As a result of this, the method overall requires 17 iterations to converge to the precision goal within a computation time of $5.1115\cdot10^{-4}$ s. The gMGF method does not show such a behavior. Instead, we observe a rapid convergence within five iteration steps which overall take $2.9384\cdot10^{-4}$ s. From Fig~\ref{Gro}.(d) we see that the the difference in computing time further increases for higher values of $y$.

The aforementioned deterioration within Newton's method shows up in the gain, which for the first iteration step amounts to $g_1\New\approx-0.78434$. After this, the gain rises to $g_2\New\approx0.10909$ from which it permanently decreases until convergence. As before, no loss is detected for the gMGF method. Instead, for the first iteration step we find $g\gMGF_1\approx0.75441$. From this the gain decreases until convergence is reached.

Further details on the convergence behavior for $y=3$ can be obtained from Tab.~\ref{tabular:Grow}. The first iteration within Newton's method is connected to an initial loss, from which a big negative value for $\varrho_2\New$ results. From the quantity $|x_n-x_{n-1}|$ we find a slight increase for $n\in[2,11]$. Within the consecutive steps, this quantity however decreases until convergence is reached. For the COC  we observe a slow convergence towards the theoretical value. In contrast to this, the gMGF approach shows a considerable magnitude in the difference between $x_1$ and $x_0$. Consequently a rapid convergence is detected also for this example. The COC varies over the five steps encountered until it tends towards the theoretical value. 

\subsection{Calculating the coolant mass flow in a heat exchanger}
Heat exchangers are used in many applications of engineering such as district heating, chemical processing or automotive thermal management. The modeling of these systems is of interest, in particular, when on-demand temperature control strategies are to be utilized in order to operate technical systems at an optimal temperature. Furthermore, if the technical system under consideration has to be monitored, a mathematical modeling of heat exchangers may also be of interest (c.f. Ref.~\cite{Herzog, Herzog17}). 

Here, we want to restrict the discussion to the case, that a heat exchanger is used in order to cool a working fluid via thermal contact with a coolant. The medium to be cooled which is transmitted with the mass flow $\dot m_1$ through the heat exchanger exhibits the temperature $T_\mathrm{1,u}$ upstream to the device. The corresponding quantities of the cooling fluid are given by $\dot m_2$ and $T_\mathrm{2,u}<T_\mathrm{1,u}$. Due to the thermal exchange, the working fluid temperature $T_\mathrm{1,d}$ downstream to the device results. Assuming single-phase fluid flow, under stationary conditions the thermal behavior of a thermally well-isolated heat exchanger may be quantified by means of its dimensionless temperature change $P$. This quantity relates the temperature difference in the working fluid to the maximum stationary temperature difference in the system, hence~\cite{VDI13}
\begin{equation}\label{DTC}
	P=\frac{T_\mathrm{1,u}-T_\mathrm{1,d}}{T_\mathrm{1,u}-T_\mathrm{2,u}}.
\end{equation}
Physically, the dimensionless temperature change may be regarded as the efficiency of the heat exchanger. Its concrete value depends on the current working point. In the context of convective thermal exchange, $P\in[0,1)$ can be expressed as a function of the heat capacity flow rates $\Gamma_\mathrm{f}=c_\mathrm{p,f}\dot m_\mathrm{f}$, where $c_\mathrm{p,f}$ is the isobaric specific heat capacity of fluid $\mathrm{f}\in\{1,2\}$, hence $P=P(\Gamma_1,\Gamma_2)$. Assuming that the temperatures are known from sensors, the objective of the present example is to resolve the identity
\begin{equation}\label{ObjectCC}
	\frac{T_\mathrm{1,u}-T_\mathrm{1,d}}{T_\mathrm{1,u}-T_\mathrm{2,u}}=P(\Gamma_1,\Gamma_2)
\end{equation}
for $\Gamma_2$, which allows to determine the cooling fluid mass flow. The solution to this problem is of interest both for control purposes, as well as for diagnostic functions, where e.g. a leakage in the thermal system may be detected, if $\Gamma_2$ calculated from Eq.~(\ref{ObjectCC}) does not fit the measured value.

\begin{figure}[t]
\begin{center}
	\includegraphics[width=0.6\linewidth]{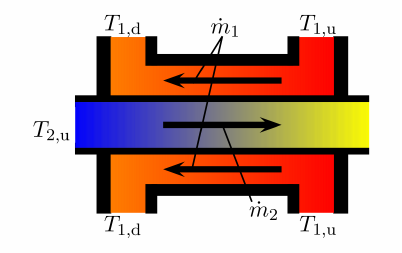}
	\caption{Schematic representation of a heat exchanger with counter current flow. The working fluid of mass flow $\dot m_1$ is characterized by the temperatures $T_\mathrm{1,u}$ and $T_\mathrm{1,d}$ up- and downstream to the device respectively. The cooling fluid exhibits the mass flow $\dot m_2$ with a coolant temperature $T_\mathrm{2,u}$ upstream to the heat exchanger.}\label{HEX}
\end{center}
\end{figure}
Apart from the heat capacity flow rates, $P$ also depends on the exact geometry of the heat exchanger. Here we will investigate a counter current flow arrangement for the mass flows $\dot m_1$ and $\dot m_2$ as indicated in Fig.~\ref{HEX} by the respective arrows. The temperatures, occurring within the dimensionless temperature change (c.f. Eq.~(\ref{DTC})) are also marked within this picture. Setting the left hand side of Eq.~(\ref{ObjectCC}) equal to $y$ and defining $x=\Gamma_1/\Gamma_2$, for such a device, we end up at
\begin{equation}\label{PCC}
	y=\frac{1-\exp\left\{-\frac{\gamma(\Gamma_1,x)}{\Gamma_1/\Gamma_{1,0}}\left(1-x\right)\right\}}{1-x\exp\left\{-\frac{\gamma(\Gamma_1,x)}{\Gamma_1/\Gamma_{1,0}}\left(1-x\right)\right\}},
\end{equation}
where we have introduced the auxiliary quantity $\Gamma_{1,0}=1$~J/(Ks) to avoid discussions on unities in the following. If the thermal resistance of the material as well as temperature effects on the fluid properties can be neglected, the total thermal conductivity of the two fluids and the heat exchanger can effectively be written as~\cite{Herzog17}
\begin{equation}\label{gammaCC}
	\gamma(\Gamma_1,x)=\frac{a\ \left(\Gamma_1/\Gamma_{1,0}\right)^\mu}{b\ x^\mu+1}.
\end{equation}
The positive numbers $a$, $b$, and $\mu$ serve as parameters which take account of the thermal behavior of the heat exchanger. For the current simulation, we used the values $a=10$, $b=12$ and $\mu=0.8$, as well as $\Gamma_1=250$ J/(Ks). We chose to invert the model for $y\in[0.10,0.99]$ with the step size $\delta y=0.01$ and the initial value $x_0=2.5$.

\begin{figure}[t]
	\begin{center}
	\includegraphics[width=.99\linewidth]{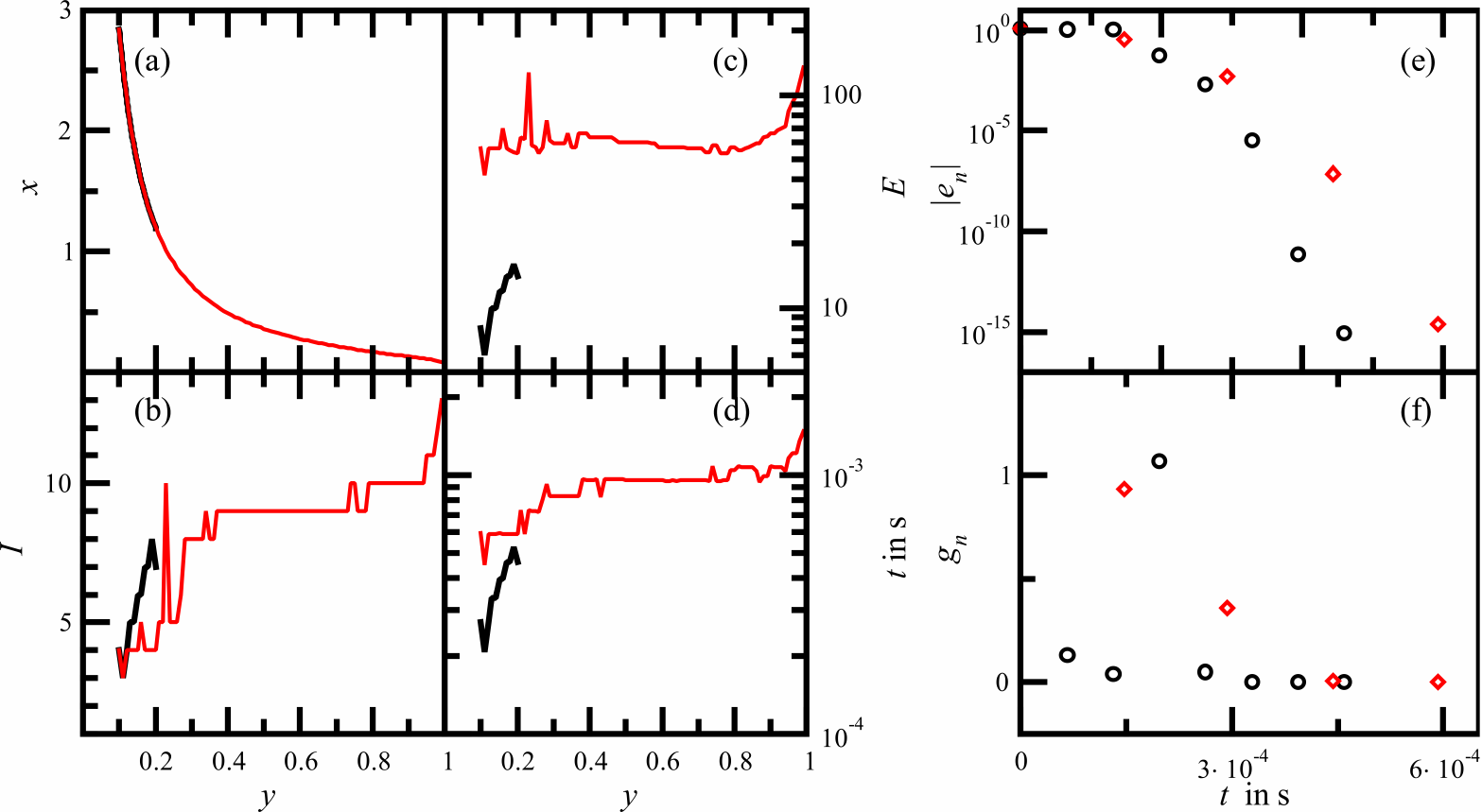}
	\caption{Comparison of results for solving Eq.~(\ref{PCC}) with Eq.~(\ref{gammaCC}) as obtained from Newton's method (black curve) and our approach (red curve). For the shown outcomes the values $a=10$, $b=12$, $\mu=0.8$, $\Gamma_{1,0}=1$~J/(Ks) and $\Gamma_1=250$ J/(Ks) have been used: subplot (a) displays the inverted function $x=x(y)$, whereas (b) shows the number of iterations needed for convergence. The number of function evaluations connected to this can be seen from subplot (c). In subplot (d) the computing time until convergence is shown. We see that Newton's method breaks down for $y>0.2$. Criteria to compare the convergence behavior of Newton's method and our approach for  solving Eq.~(\ref{PCC}) with $y=0.2$ are shown in subplots~(e) and~(f). The circles display results as obtained from Newton's method, whereas the diamonds are calculated from the approach suggested in this publication. In subplot (e) the magnitude of the residuals as a function of the computing time is shown. The gain is displayed in subplot (f).}\label{Hea}
	\end{center}
\end{figure}

Results from these investigations can be taken from Fig.~\ref{Hea}. As we see from subplot~(a), contrary to the previous application examples, only the iterative gMGF approach succeeds in covering the full interval for $y$ taken into account. In contrast to this, Newton's method exhibits convergence only in the range $[0.1,0.2]$. A closer investigation reveals, that for $y>0.2$ we have $x_1\New<0$. As a consequence of this $x_n\New$ with $n\geq2$ turn out to be complex, where both, the real- as well as the imaginary part diverge during the iteration. We suspect, that Newton's method is unable to handle the strong nonlinear behavior of the problem at hand. Due to its instability, it could not be used in the industrial application at hand. To the contrary no critical behavior is observed within the gMGF method.

The number of iterations displayed in subplot~(b) reveal that $I\New(y)=I\gMGF(y)$ for $y\in\{0.1,0.11,0.12\}$. For all other values the gMGF approach outperforms Newton's method as far as the number of iterations are concerned.

From subplot~(c), we see, that Newton's method clearly outperforms the gMGF approach in the region, where both algorithms converge. Moreover, by comparing $I\gMGF$ to $E\gMGF$, we may deduce that the magnitude of the degree in the considered range is quite high. A detailed analysis reveals, that in the considered interval $\varkappa(x)$ ranges between 1 and -22. Bearing in mind, that  $\mathfrak{H}_\varkappa$ is utilized to reformulate the initial problem by means of an essentially linear power series, the conjecture that the failure of Newton's method is due to the strong nonlinear nature of the function under consideration is strongly supported.

Turning towards subplot~(d), we find, as with the number of function evaluations, that Newton's method, if it converges, is superior to the iterative gMGF approach. For the latter we observe a tendency of $t\gMGF(y)$ which resembles the behavior of $E\gMGF(y)$, even though it is less pronounced. This again is due to the fact, that the operations of determining $\varkappa$ and evaluating $\mathfrak{H}_{\varkappa}$ from $\mathfrak{H}_0$ come with a low computational cost.

For the investigation of of the magnitude of the error displayed in subplot~(e), we have chosen $y=0.2$. We see, that both methods converge rapidly with an advantage on the side of Newton's method. There is no deterioration, as was observed within the previous applications. Finally turning towards the gain, we see, that this value reaches its maximum value in the third iteration step for Newton's method, whereas within the gMGF algorithm, the first step yields the highest figure.
\begin{table}[t]
	\centering
	\caption{Numerical results for solving Eq~(\ref{PCC}) with Eq.~(\ref{gammaCC}) for $x$ with $y=0.2$ and $x_0=2.5$}
	\label{tabular:T_Hot}
	\begin{tabular}[!ht]{c c c c c c c c c c c c c c c}
	\toprule[2pt]
		Method\\
		\midrule[1.5pt]
		Newton & $n=1$ & $n=2$ & $n=3$ & $n=4$ & $n=5$ & $n=6$  & $n=7$ & $I$\\
		\midrule
		$|x_n-x_{n-1}|$ & $2.5$ & $0.040$ & $1.1$ & $0.052$ & $2.1\cdot10^{-3}$ & $3.2\cdot10^{-6}$ & $7.3\cdot10^{-12}$ & 7 \\
		$\varrho_n$ &  & $.321$ & $87.6$ & $1.07$ & $2.00$ & $2.00$ & $2.00$ \\\midrule[1.5pt]
		&   &  &  &  &  & &  & CPU time & &\\\midrule
		 &  &  &  &  &  & & & $4.02\cdot10^{-4}$ s\\\midrule[1.5pt]
		gMGF & $n=1$ & $n=2$ & $n=3$ & $n=4$ &  &  & $I$ & CPU time\\
		\midrule
		$|x_n-x_{n-1}|$\ \  & $1.7$ & $0.36$ & $5.0\cdot10^{-3}$ & $7.2\cdot10^{-8}$ &  &  & $4$ & $4.74\cdot 10^{-4}$ s\\
		$\varrho_n$ &  & $3.35$ & $2.61$ & $2.98$ &  & \\
		\bottomrule[2pt]
	\end{tabular} 
\end{table}

Finally in Tab.~\ref{tabular:T_Hot}, we can obtain further information on the convergence of the method for $y=0.2$. Newton's method does not exhibit any losses at any stage of the iteration. This is due to the fact, that the actual solution is sufficiently close to the initial value. From $|x_n-x_{n-1}|$ we see that Newton's method converges rapidly within seven steps. The COC exhibits a huge value which is connected to the big difference between $|x_3-x_2|$ and $|x_2-x_1|$. We also observe, that $\varrho_2\New$ tends towards the corresponding theoretical value. In contrast within the gMGF method overall four iterations are required to reach convergence. Due to the higher computation time connected to one iteration step, as far as CPU time is concerned Newton's method outperforms the gMGF algorithm in this example. Interestingly, for $\varrho_n\gMGF$ we do not observe a trend towards 2.00. Instead, we find that $\varrho_4\gMGF=2.98$. This is due to the circumstance, that the expression $|h_{0,2}(x_n)/h_{1,2}^2(x_n)|\approx10.9934$ for $n\in\{2,3,4\}$. From this, we conclude that $h_{\kappa,2}(x_n)$ is very small, such that the the expression on the right hand side of Eq.~(\ref{errorEq}) turns out to be smaller when compared to higher order terms. We have checked, that with further iterations aiming at a more demanding precision goal, the COC of the gMGF finally tends towards the theoretical value.

\section{Summary and outlook}\label{SumOut}
Summarizing, in this paper, we have developed a methodology to solve nonlinear equations where the computing time has a low sensitivity on the initial value. For many applications in engineering, this may be a crucial asset, in particular if safety issues are to be taken into consideration. Instead of increasing the local computational efficiency, we chose a novel approach by introducing the concept of generalized moment generating functions. From this we may systematically transform the problem to a form which is essentially represented by a linear approximation. As a result of this, the global convergence is increased such that a small number of iterations is required to reach the precision goal. In this publication, we have exemplified the approach by generalizing Newton's method.

We have demonstrated the methodology at the hand of seven numerical examples and three applications. Especially for high deviations $|\zeta-x_0|$ between the initial value $x_0$ and the actual solution $\zeta$, the generalized moment generating function approach outperforms Newton's method. For the considered examples and applications we found an essentially flat computing time within the proposed method hinting towards an improvement of the global convergence with respect to Newton's method. In particular, we observed, that a reduced precision goal favors the performance of the presented algorithm.

Both these aspects, the low sensitivity of the computing time with respect to the starting value, as well as the rapid convergence to a crude precision target, are requirements which are met in a number of applications, such as real-time modeling of power systems. Other examples include diagnosis or control in process or automotive engineering. Here, the corresponding state variables like temperatures, pressures and mass flows are to be determined with a strict real-time requirement and a precision target of usually only a few significant digits. Especially in the case of mobile microcontrollers only very restricted computing powers are available. In this respect, the discussed approach could indeed enable the implementation of novel numerical solutions which have not been computationally affordable so far.

We are confident, that the discussed approach could also be applied to other iteration algorithms like those, which have been mentioned in the introduction. This offers the prospect to also speed up these iterations and thus decrease computing times for equation solving in general. A concrete generalization of the respective algorithms requires, of course, further detailed investigations. So one obvious direction for future research consists in trying to generalize existing one- as well as multi-step methods within the framework of generalized moment generating functions. This should also be possible for multidimensional problems.

Moreover, it may be possible to apply the discussed method in the framework of control problems. We think, that the rapid convergence of the method may be beneficial especially in the field of automated parametrization of controllers or as solver in the framework of model based control. One concrete example is to use the method in order to speed up convergence in reinforcement learning of model predictive control. As discussed in \cite{BERTSEKAS2022100121} a common strategy in this directions is to reformulate the occurring one-step look-ahead minimization problem as the solution of Bellmann's equation. The generic solution to this through Newton's method requires a proper initial guess. As we have seen in this publication, the iteration based on the generalized moment generating function is far more robust in this respect. Therefore, a possibility to improve reinforcement learning could consist in applying the proposed method in this context.

\end{document}